\documentclass[twocolumn]{revtex4}
\usepackage{epsf}
\usepackage{url}
\usepackage{graphicx}
\usepackage{boxedminipage}
\usepackage{bm}
\usepackage{amsmath, amsthm, amssymb} 

\def\coma#1{{\bf [[ Aleks: #1 ]]}}

\def\comaOff#1{}
\def\comaWOff#1{}
\def\tfidf{{\tt tf*idf}}
\def\mat#1{\vec{#1}}

\newtheorem{lemma}{Lemma}

\def\vec#1{\ensuremath{\mathchoice
                     {\mbox{\boldmath$\displaystyle#1$}}
                     {\mbox{\boldmath$\textstyle#1$}}
                     {\mbox{\boldmath$\scriptstyle#1$}}
                     {\mbox{\boldmath$\scriptscriptstyle#1$}}}}%
\def\tens#1{\relax\ifmmode\mathsf{#1}\else\textsf{#1}\fi}

\begin{document}

\title{Discrete Component Analysis}

\author{Wray Buntine}
\affiliation{Helsinki Institute for Information Technology (HIIT).
Dept.\ of Computer Science, PL 68. 00014, University of Helsinki,
Finland.} \email{Wray.Buntine@hiit.fi}
\author{Aleks Jakulin}
\affiliation{Columbia University, Department of Statistics, 1255
Amsterdam Avenue, New York, NY 10027-5904, USA.}
\email{jakulin@acm.org}

\begin{abstract}%
This article presents a unified theory for analysis of components
in discrete data, and compares the  methods with techniques such
as independent component analysis,  non-negative matrix
factorisation and  latent Dirichlet allocation.  The
main families of algorithms discussed are a variational approximation, Gibbs
sampling, and Rao-Blackwellised Gibbs sampling. Applications are
presented for voting records from the United States Senate for
2003, and for the Reuters-21578 newswire collection.
\end{abstract}

\maketitle

\section{Introduction}

Principal component analysis (PCA) \cite{Mardia.book}
is a key method in the statistical engineering toolbox. It is well over a
century old, and is used in many different ways.
PCA is also known as the Karh\"{u}nen-Lo\`{e}ve transform
or Hotelling transform in image analysis,
and a variation is latent semantic analysis (LSA) in text analysis
\cite{deerwester90indexing}.
It is a kind
of eigen-analysis since it manipulates the eigen-spectrum of the
data matrix.
It is usually applied to
measurements and real valued data, and used for feature extraction
or data summarization.
LSA might not perform the centering step
(subtracting the mean from each data vector prior to eigen-analysis)
on the word counts
for a document
to preserve matrix sparseness, or might convert the word counts
to real-valued \tfidf\ \cite{baezayates99}.
The general approach here is {\it data reduction}.

Independent component analysis (ICA, see \cite{ICA}) is in some ways
an extension of this general approach, however it also involves the
estimation of so-called latent, unobservable variables. This kind of
estimation follows the major statistical methodology that deals with
general unsupervised methods such as clustering and factor analysis.
The general approach is called {\it latent structure analysis},
which is more recent, perhaps half a century old. The data is
modelled in a way that admits unobservable variables, that influence
the observable variables. Statistical inference is used to
``reconstruct'' the unobservable variables from the data jointly
with general characteristics of the unobservable variables
themselves. This is a theory with particular assumptions (i.e., a
``model''), so the  method may arrive at poor results.

Relatively recently the statistical computing and machine learning
community has become
aware of seemingly similar approaches for discrete observed data that
appears under many names.  The best known of
these in this community are
probabilistic latent semantic indexing (PLSI)
\cite{hofmann99probabilistic}, non-negative matrix factorisation
(NMF) \cite{lee99learning} and  latent Dirichlet allocation (LDA)
\cite{Blei03}.  Other variations are discussed later in
Section~\ref{sct-hn}.
We refer to these methods jointly as {\em Discrete Component Analysis} (DCA),
and this article provides a unifying model for them.

All the above
approaches assume that the data is formed from individual
observations (documents, individuals, images), where each
observation is described through a number of variables (words,
genes, pixels). All these approaches attempt to summarize or
explain the similarities between observations and the correlations
between variables by inferring latent variables for each
observation, and associating latent variables with observed
variables.

These methods are applied in the social sciences, demographics and
medical informatics, genotype inference, text and image analysis, and
information retrieval. By far the largest body of applied work in this
area (using citation indexes) is in genotype inference due to the
Structure program \cite{pritchard:1999}. A growing body of work is in
text classification and topic modelling (see
\cite{griffiths04pnas,buntine2004sawm}), and language modelling in
information retrieval (see
\cite{Azzopardi2003,BuntineJakulin04,canny04}).
As a guide, argued in the next section,
the methods apply when PCA or ICA might be used, but the
data is discrete.

Here we present in Section~\ref{sct-dfn} a unified theory for analysis
of components in discrete data, and compare the methods with related
techniques in Section~\ref{sct-hn}.  The main families of algorithms
discussed in Section~\ref{sct-alg} are a variational approximation,
Gibbs sampling, and Rao-Blackwellised Gibbs sampling. Applications are
presented in Section~\ref{sct-app} for voting records from the United
States Senate for 2003, and the use of components in subsequent
classification.

\section{Views of DCA}

One interpretation of the DCA methods is that they are a way of
approximating large sparse discrete matrices. Suppose we have a
$500,000$ documents made up of $1,500,000$ different words. A
document such as a page out of Dr.\ Seuss's {\it The Cat in The
Hat}, is first given as a {\it sequence of words}.
\begin{quote}
So, as fast as I could, I went after my net.  And I said, ``With my net
I can bet them I bet, I bet, with my net, I can get those Things yet!''
\end{quote}
It can be put in the {\it bag of words} representation, where word
order is lost.  This yields a list of words and their counts in brackets:
\begin{quote}
after(1) and(1) as(2) bet(3) can(2) could(1) fast(1) get(1) I(7)
my(3) net(3) said(1) so(1) them(1) things(1) those(1) went(1)
with(2) yet(1) ~.
\end{quote}
Although the word `you' never appears in the original,
we do not include `you (0)' in the representation since zeros are suppressed.
This sparse vector can be represented as a vector in full word
space with $1,499,981$ zeroes and the counts above making the
non-zero entries in the appropriate places. Given a matrix made up
of rows of such vectors of non-negative integers dominated by
zeros, it is called here a {\em large sparse discrete matrix}.

Bag of words is a basic representation in information retrieval
\cite{baezayates99}. The alternative is a sequence of words. In
DCA, either representation can be used and the models act the
same, up to any word order effects introduced by incremental
algorithms.  This detail is made precise in subsequent sections.

In this section, we argue from various perspectives that large
sparse discrete data is not well suited to standard PCA or ICA methods.

\subsection{Issues with PCA}

PCA has been normally applied to numerical
data, where individual instances are vectors of real numbers.
However, many practical datasets are based on vectors of integers,
non-negative counts or binary values. For example, a particular
word cannot have a negative number of appearances in a document.
The vote of a senator can only take three values: Yea, Nay or Not
Voting. We can transform all these variables into real numbers
using \tfidf, but this is a linear weighting that
does not affect the shape of a distribution.

With respect to modelling count data in
linguistic applications, Dunning makes the following
warning \cite{dunning94accurate}:
\begin{quote}
Statistics based on the assumption of normal distribution
are invalid in most cases of statistical text analysis unless
either enormous corpora are used, or the analysis is restricted to
only the very most common words (that is, the ones least likely to
be of interest). This fact is typically ignored in much of the
work in this field. Using such invalid methods may seriously
overestimate the significance of relatively rare events.
Parametric statistical analysis based on the binomial or
multinomial distribution extends the applicability of statistical
methods to much smaller texts than models using normal
distributions and shows good promise in early applications of the
method.
\end{quote}
While PCA is not always considered a method based on Gaussians,
it can be justified using Gaussian distributions
\cite{roweis98em,tipping-bishop-99}.
Moreover, PCA is justified using a least squares distance
measure, and most of the properties of Gaussians
follow from the distance measure alone.
Rare events correspond to points far away under an $L_2$ norm.

\begin{figure*}
\begin{center}
\input{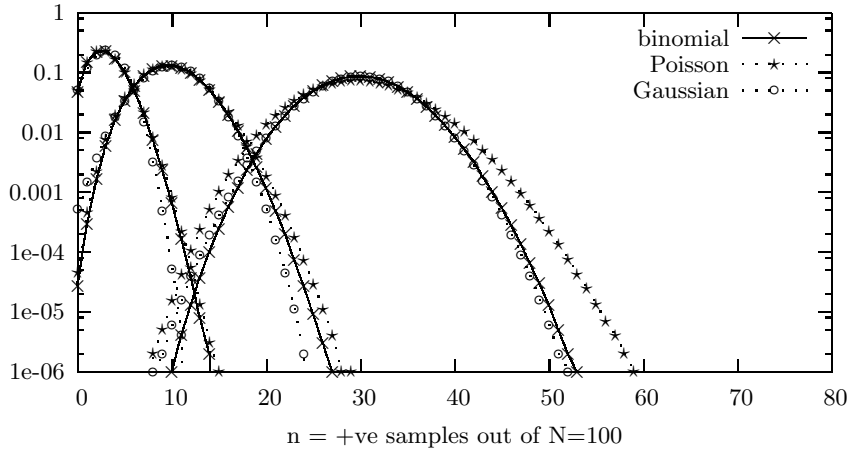}
\end{center}
\caption{Gaussian and Poisson approximations to a
binomial.\label{fig-pg100}}
\end{figure*}

Fundamentally, there are two different kinds of large sample
approximating distributions that dominate discrete statistics: the
Poisson and the Gaussian.  For instance, a large sample binomial is
approximated as a Poisson \footnote{This is a distribution on
integers where a rate is given for events to occur, and the
distribution is over the total number of events counted.} when the
probability is small and as a Gaussian otherwise \cite{Ross.intro}.
Figure~\ref{fig-pg100} illustrates this by showing the Gaussian and
Poisson approximations to a binomial with sample size $N=100$ for
different proportions ($p=0.03, 0.01, 0.03$).   Plots are done with
probability in log scale so the errors for low probability values
are highlighted.
One can clearly see the problem here:
the Gaussian provides a reasonable approximate for medium
values of the proportion $p$ but for small values it severely
underestimates low probabilities.  When these low probability events
occur, as they always will, the model becomes distorted.

Thus in image analysis
based on analogue to digital converters, where data is counts,
Gaussian errors can sometimes be assumed, but the Poisson should
be used if counts are small.  DCA then avoids Gaussian modelling
of the data, using a Poisson or multinomial directly.

Another critique of the general style of PCA comes from the
psychology literature, this time it is used as a justification for
DCA \cite{Griffiths02}. Griffiths and Steyvers argue against the
least squares distance of PCA:
\begin{quote}
While the methods behind LSA were novel in scale and subject,
the suggestion that similarity relates to distance in
psychological space has a long history (Shepard, 1957). Critics
have argued that human similarity judgments do not satisfy the
properties of Euclidean distances, such as symmetry or the triangle
inequality. Tversky and Hutchinson (1986) pointed out that
Euclidean geometry places strong constraints on the number of
points to which a particular point can be the nearest neighbor,
and that many sets of stimuli violate these constraints.
\end{quote}
They also considered power law arguments which PCA violates
for associated words.

\comaOff{PCA may yield the ``horseshoe'' effect. This is explained in
\url{http://ordination.okstate.edu/PCA.htm}. Perhaps we could put
in here an example of the horseshoe effect that would persuade
machine learning folk as realistic? Any ideas? In addition to
this, `rotation of factors' (such as varimax rotation, quartimax
rotation) is normally used in statistics to add more meaning to
components.\\ The LSI people admit that the dimensions are not
interpretable --- cite Deerwester90. - but that's not their
objective, they just want to reduce dimensionality and desparsify
the data. }
\comaWOff{OK, didn't fit this in.}


\subsection{Component Analysis as Approximation}
\label{sct-va}

In the data reduction approach for PCA, one seeks to
reduce each $J$-dimensional data vector to a smaller $K$-dimensional vector.
This can be done by approximating the full data matrix as a
product of smaller matrices, one representing the
reduced vectors called the component/factor {\it score matrix},
and one representing a data independent part called
the component/factor  {\it loading matrix},
as shown in Figure~\ref{fig-mtx}.
In PCA according to least squares theory,
this approximation is made by eliminating the lower-order
eigenvectors, the least contributing components \cite{Mardia.book}.
\begin{figure*}
\[
\begin{array}{cccccccc}
 & {J \textrm{ words}} & & {K \textrm{ components}} & & &  & \\
 \rotatebox{90}{\hspace{-1cm}I \textrm{ documents}} &
\hspace{6mm} \overbrace{\hspace{-6mm}\left\{
\left(\begin{array}{cccccc}
w_{1,1} & w_{1,2} & & \cdots & & w_{1,J}\\
w_{2,1} & w_{2,2} & & \cdots & & w_{2,J}\\
\\
\vdots & \vdots & & \ddots & & \vdots\\
\\
w_{I,1} & w_{I,2} & & \cdots & & w_{I,J}\\
\end{array}\right)\right.\hspace{-4mm}} \hspace{4mm}
 & \simeq &
\hspace{3mm} \overbrace{\hspace{-3mm} \left(\begin{array}{ccc}
l_{1,1} & \cdots & l_{1,K}\\
l_{2,1} & \cdots & l_{2,K}\\
\\
\vdots & \ddots & \vdots\\
\\
l_{I,1} & \cdots & l_{I,K}\\
\end{array}\right)\hspace{-4mm}} \hspace{4mm}
& * & \begin{array}{c}{J \textrm{ words}}\\ \hspace{3mm}
\overbrace{\hspace{-3mm}\left. \left(\begin{array}{cccccc}
\theta_{1,1} & \theta_{2,1} & & \cdots & & \theta_{J,1}\\
\vdots & \vdots   & & \ddots & & \vdots\\
\theta_{1,K} & \theta_{2,K} & & \cdots & & \theta_{J,K}\\
\end{array}\right)\right\}\hspace{-6mm}} \hspace{6mm}\\ ~ \end{array}
& \rotatebox{90}{\hspace{-1.2cm}K \textrm{ components}}
\\ \\
& \textrm{data matrix} & & \textrm{score matrix} & &
\textrm{loading matrix}^T &
 \end{array}  \]
\caption{The matrix approximation view} \label{fig-mtx}
\end{figure*}

If there are $I$ documents, $J$ words and $K$ components, then the
matrix on the left has $I*J$ entries and the two matrices on the
right have $(I+J)*K$ entries. This represents a simplification
when $K\ll I,J$. We can view DCA methods as seeking the same goal
in the case where the matrices are sparse and discrete.

When applying PCA to large sparse discrete matrices, or LSA using word
count data interpretation of the components, if
it is desired, becomes difficult
(it was not a goal of the original method \cite{deerwester90indexing}).
Negative values appear in the
component matrices, so they cannot be interpreted as ``typical
documents'' in any usual sense. This applies to many other kinds of
sparse discrete data: low intensity images (such as astronomical
images) and verb-noun data used in language models introduced by
\cite{pereira1993}, for instance.

The cost function being minimized then plays an important role.
DCA places constraints on the
approximating score matrix and loading matrix in Figure~\ref{fig-mtx}
so that they are also non-negative.
It also uses an entropy distance instead of a least squares
distance.

\subsection{Independent Components}
\label{sct-ica}

Independent component analysis (ICA) was also developed as an
alternative to PCA.
Hyv\"{a}nen and Oja \cite{hyvaoja2000} argue
that PCA methods merely find uncorrelated
components.
ICA then was
developed as a way of representing multivariate data with truly
{\it independent} components.
In theory, PCA approximates this also if the data is Gaussian
\cite{tipping-bishop-99}, but in practice it rarely is.

The basic formulation is that a $K$-dimensional data vector
$\vec{w}$ is a linear invertible function of $K$ independent
components represented as a $K$-dimensional latent vector $\vec{l}$,
$\vec{w} \,=\,{\mat\Theta} \vec{l}$ for a square invertible matrix
$\mat\Theta$. Note the ICA assumes $J=K$ in our notation.
$\mat\Theta$  plays the same role as the loading matrix above. For
some univariate density model U, the independent components are
distributed as $p(\vec{l}\,|\, U)  \,=\, \prod_k p(l_k\,|\, U)$,
thus one can get a likelihood formula $p(\vec{w}\,|\, {\mat\Theta},U
)$ using the above equality \footnote{By a change of coordinates
\[
    p(\vec{w}\,|\, {\mat\Theta},U ) ~=~ \frac{1}{\det(  {\mat\Theta} )}
   \prod_k p\left( \left( {\mat\Theta}^{-1} \vec{w} \right)_k\,|\, U \right)
\]}.

The Fast ICA algorithm  \cite{hyvaoja2000} can be interpreted as
a maximum likelihood approach based on this model and likelihood formula.
In the sparse discrete case, however, this formulation breaks down
for the simple reason that $\vec{w}$ is mostly zeros: the equation
can only hold if $\vec{l}$ and ${\mat\Theta}$ are discrete as well
and thus the gradient-based algorithms for ICA cannot be
justified. To get around this in practice, when applying ICA to
documents \cite{bingham03}, word counts are sometimes first turned
into \tfidf\ scores \cite{baezayates99}.

To arrive at a formulation more suited to discrete data, we can
relax the equality in ICA (i.e., $\vec{w} \,=\,{\mat\Theta}
\vec{l}$) to be an expectation:
\[
\mathbb{E}_{\vec{w} \sim p(\vec{w}|\vec{l},U)}\left[
 \vec{w}\right] ~ = ~ {\mat\Theta}  \vec{l}~.
\]
We still have independent components, but a more robust
relationship between the data and the score vector. Correspondence between ICA
and DCA has been noted in \cite{BuntineJakulin04,canny04}.
With this expectation relationship, the dimension of
$\vec{l}$ can now be less than  the dimension of $\vec{w}$,
$K<J$, and
thus ${\mat\Theta}$ would be a rectangular matrix.

\section{The Basic Model}
\label{sct-dfn}

A good introduction to these models from a number of viewpoints is
by \cite{Blei03,canny04,BuntineJakulin04}. Here we present a
general model. The
notation of words, bags and documents will be used throughout,
even though other kinds of data representations also apply.
In statistical terminology, a word is an observed variable,
and a document is a data vector (a list of observed variables)
representing an instance.
In machine learning terminology, a word is a feature, a bag is a data
vector, and a document is an instance.  Notice that the bag collects
the words in the document and loses their ordering. The bag is
represented as a data vector $\vec{w}$. It is now $J$-dimensional.
The latent, hidden or unobserved vector $\vec{l}$ called the component {\it scores} is
$K$-dimensional. The term {\it component} is used
here instead of topic, factor or cluster. The parameter matrix is the
previously mentioned component loading matrix ${\mat\Theta}$,
and  is $J \times K$.

At this point, it is also convenient to introduce
the symbology used throughout the paper.
The symbols summarised in Table~\ref{tbl-syms} will be introduced
as we go.
\begin{table*}
\begin{center}
\begin{tabular}{l|ll}
$I$   & & number of documents\\
$(i)$   & & subscript to indicate document, sometimes dropped\\
$J$   & & number of different words, size of the dictionary\\
$K$   & & number of components\\
$L_{(i)}$ & & number of words in document $i$\\
$S$   & & number of words in the collection, $\sum_i L_{(i)}$ \\
$\vec{w}_{(i)}$ & & vector of $J$ word counts in document $i$, row totals of
       $\mat V$, entries $w_{j,(i)}$ \\
$\vec{c}_{(i)}$ & & vector of $K$ component counts for document $i$,
column  totals of  $\mat V$\\
$\mat V$ & & matrix of word counts per component, dimension
   $J\times K$, entries $v_{j,k}$\\
$\vec{l}_{(i)}$ & & vector of $K$ component scores for document $i$, entries $l_{k,(i)}$ \\
$\vec{m}_{(i)}$ & & $\vec{l}_{(i)}$ normalised, entries $m_{k,(i)}$ \\
$\vec{k}_{(i)}$ & & vector of $L_{(i)}$ sequential
   component assignments for the words in \\
&& document $i$,
entries $k_{l,(i)}\in [1,\ldots,K]$\\
$\mat\Theta$ & & component loading matrix, dimension $J\times K$, entries
       $\theta_{j,k}$\\
$\vec\theta_{\cdot,k}$ & & component loading vector for component $k$,
       a column of $\mat\Theta$ \\
$\vec\alpha,\vec\beta$ & & $K$-dimensional parameter vectors for
      component priors\\
\end{tabular}
\end{center}
\caption{Summary of major symbols}
\label{tbl-syms}
\end{table*}

\comaOff{Second, you use two different symbols for score/loading
depending on whether the model is GP or MPCA. I don't think we
need two different symbols here, but even if we choose to have
them, it's important to explain the similarity.}
\comaWOff{The $\vec{l}$ is unnormalised, and the $\vec{m}$ normalised.
Would be good to do this, but would take time.
I'll get to it later if time permits.}
\comaOff{how about referring to a document with $\vec{x}$
instead of $\vec{d}$? How about referring to the latent
description of a document with $\vec{z}$ instead of $\vec{l}$?}
\comaWOff{trouble is that depending on the formulation, you can add different
hidden variables, e.g., $\mat{V}$ as well, so couldn't justify
neatly using $z$, and thus $y$ got dropped too}

\subsection{Bags or Sequences of Words?}
\label{sct-bs}

For a document $\vec{x}$ represented as a sequence of words, if
$\vec{w}=\mbox{bag}(\vec{x})$ is its bagged form, the bag of words,
represented as a vector of counts. In the simplest case, one can use
a multinomial with sample size $L=|\vec{x}|$ and vocabulary size
$J=|\vec{w}|$ to model the bag, or alternatively $L$ independent
discrete distributions \footnote{The {\it discrete} distribution is
the multivariate form of a Bernoulli where an index
$j\in\{0,1,\ldots,J-1\}$ is sampled according to a $J$-dimensional
probability vector.} with $J$ outcomes to model each $x_l$. The bag
$\vec{w}$ corresponds to the sequence $\vec{x}$ with the order lost,
thus there are $\frac{\left(\sum_j w_j\right)!}{\prod_j w_j ! }$
different sequences that map to the same bag $\vec{w}$. The
likelihoods for these two simple models thus differ by just this
combinatoric term.

Note that some likelihood based methods such as maximum
likelihood, some Bayesian methods, and some other fitting methods
(for instance, a cross validation technique) use the likelihood as
a black-box function.  They take values or derivatives but otherwise do
not further interact with the likelihood. The combinatoric term
mapping bag to sequence representations can be
ignored here safely because it does not affect the fitting of the
parameters for $\cal M$. Thus for these methods, it is irrelevant
whether the data is treated as a bag or as a sequence. This is a
general property of multinomial data.

\comaOff{What is a `functional'? Maybe `criterion', `loss' or
`score' would be more meaningful for machine learners.} Thus, while
we consider bag of words in this article, most of the theory applies
equally to the sequence of words representation \footnote{Some
advanced fitting methods such as Gibbs sampling do not treat the
likelihood as a black-box. They introduce latent variables that
expands the functional form of the likelihood, and they may update
parts of a document in turn. For these, ordering effects can be
incurred by bagging a document, since updates for different parts of
the data will now be done in a different order. But the combinatoric
term mapping bag to sequence representations will still be ignored
and the algorithms are effectively the same up to the ordering
affects.}. Implementation can easily address both cases with little
change to the algorithms, just to the data handling routines.

\comaOff{again, too intricate, too much detail at this point.}
\comaWOff{Made it a footnote.}

\subsection{General DCA}
\label{sct-gdca}

The general formulation introduced in Section~\ref{sct-ica} is an
unsupervised version of a linear model, and it applies to the bag
of words $\vec{w}$ as
\begin{equation}
\label{eq-dfn}
\mathbb{E}_{\vec{w} \sim p(\vec{w}|\vec{l},\mat\Theta)}\left[ \vec{w}\right] \, = \, {\mat\Theta}  \vec{l}
\end{equation}
The expected value (or mean) of the data is given by
the dot product of the component loading matrix ${\mat\Theta}$ and some latent
component scores $\vec{l}$.

\comaOff{why is it a good idea to have a likelihood model?}
In full probability (or Bayesian) modelling
\cite{Gelman.Carlin.95},
we are required to give a distribution
for all the non-deterministic values in a model, including model parameters
and the latent variables.
In likelihood modelling
\cite{casella.berger:90},
we are required to give a distribution
for all the data values in a model, including observed and
latent variables.
These are the core methodologies in computational statistics,
and most others extend these two.
The distribution for the data is called a likelihood in both
methodologies.

The likelihood of a document is the
primary way of evaluating a probabilistic model. Although
likelihood is not strictly a probability in classical statistics,
we can interpret them as
a probability that a probabilistic model $\cal M$ would generate a
document $\vec{x}$, $P(\vec{x}|{\cal M})$. On the other hand, it
is also a way of determining whether the document is usual or
unusual: documents with low likelihood are often considered to be
outliers or anomalies. If we trust our documents, low likelihoods
indicate problems with the model. If we trust out model, a low
likelihood indicates problems with a document.

Thus to complete the above formulation for DCA, we need to give
distributions matching the constraint in Equation~(\ref{eq-dfn}),
to specify the likelihood.
Distributions are needed for:
\begin{itemize}
\item
how the sequence
$\vec{x}$ or bag $\vec{w}$ is distributed
given its mean ${\mat\Theta}\vec{l}$
formed from the component loading matrix,
\item
how the component scores $\vec{l}$
are distributed,
\item
and if full probability modelling is used,
how the component loading matrix ${\mat\Theta}$
is distributed apriori,
as well as any parameters.
\end{itemize}

The formulation of Equation~(\ref{eq-dfn}) is also called an {\it
admixture model} in the statistical literature \cite{pritchard:1999}.
This is in contrast with a {\em mixture model}
\cite{Gelman.Carlin.95} which uses a related
constraint
\[
\mathbb{E}_{\vec{w} \sim p(\vec{w}|\vec{l})}\left[ \vec{w}\right] \, = \,
   \vec\theta_{\cdot,k}~,
\]
for some latent variable $k$ representing the single latent
component for $\vec{w}$.
Since $k$ is unobserved, this also corresponds to
making a weighted sum of the probability distributions
for each $\vec\theta_{\cdot,k}$.
\comaOff{the above paragraph is interesting, but again, can it be
explained in simper terms?}

\section{The Model Families}

This section introduces some forms of DCA using specific
distributions for the sequence
$\vec{x}$ or bag $\vec{w}$ and the component scores $\vec{l}$.
The fundamental model here is the Gamma-Poisson Model
(GP model for short). Other
models can be presented as variations.
The probability for a document is given for each model,
both for the case where the latent variables are
known (and thus are on the right-hand side),
and for the case where the latent variables
are included in the left-hand side.

\subsection{The Gamma-Poisson Model}
\label{sct-gp}

The general Gamma-Poisson form of DCA, introduced as GaP
\cite{canny04} is now considered in more detail:
\begin{itemize}
 \item Document data is supplied in the form of {\em word counts}.
The word count for each word type is $w_j$. Let $L$ be the total
count, so $L=\sum_j w_j$.
\comaOff{Nice!! That's good approach to explanation. But $l_k$ is
ugly - especially as $l$ is actually a vector. How about
$\vec{l}_k$ instead?}

 \item The document also has {\em component scores}
$\vec{l}$ that indicate the amount of the component in the document.
These are latent or unobserved.
The entries $l_k$ are independent and gamma
distributed
\[
l_k ~\sim~\mbox{Gamma}(\alpha_k,\beta_k) ~~~~~~~~~~~\mbox{for $k=1,\ldots,K$}.
\]
The $\beta_k$ affects scaling of the components
\footnote{Conventions for the gamma vary.  Sometimes a parameter
$1/\beta_k$ is used. Our convention is revealed in
Equation~(\ref{eq-lf}).}, while $\alpha_k$ changes the shape of the
distribution, shown in Figure~\ref{fig-gamma}.
\begin{figure*}
\begin{center}
\input{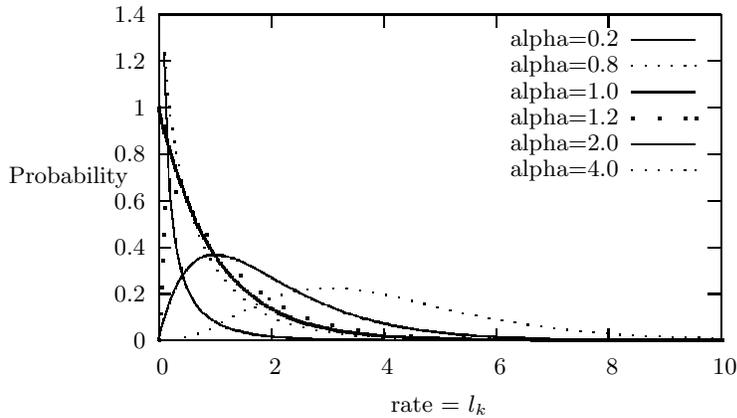}
\caption{Gamma distribution for different values of $\alpha_k$.}
\label{fig-gamma}
\end{center}
\end{figure*}
\comaOff{at this point, gamma needs to be explained. why gamma?
evidence for gamma? how does gamma look like?}

\item There is a {\em
component loading matrix} ${\mat\Theta}$ of size $J\times K$ with
entries $\theta_{j,k}$ that controls the partition of features
amongst each component. In the matrix, each column for component $k$ is
normalised across the features, meaning that $\sum_j
\theta_{j,k}=1$.
Thus each column represents the proportions of words/features in
component $k$.
\comaOff{what's th
meaning of $\alpha,\beta$? Now you're speaking of a 'general model
matrix' which was 'parameter matrix' in the prev. sections. How
about 'component loading matrix'? Normalised across the features
== for each component $k$, words are expressed as proportions, and
all the words sum up to 1.}
\item The observed data $\vec{w}$ is now Poisson distributed, for each $j$
\[
w_{j} ~\sim~ \mbox{Poisson}\left( ({\mat\Theta} \vec{l})_j \right)~.
\]
\item The $K$ parameters $(\alpha_k,\beta_k)$ to the gamma distributions give
$K$-dimensional parameter vectors
$\vec\alpha$ and $\vec\beta$.
Initially these vectors will be treated as constants, and
their estimation is the subject of later work.
\item
When using Bayesian of full probability modelling,
a prior is needed
for $\mat\Theta$. A Dirichlet prior can be used for each $k$-th
component of $\mat\Theta$ with $J$ prior parameters $\gamma_j$,
so $\vec\theta_{\cdot,k} \,\sim\, \mbox{Dirichlet}_J(\vec\gamma)$.
In practice we use a Jeffreys' prior, which has $\gamma_j=0.5$.
The use of a Dirichlet has no strong justification other than being
conjugate \cite{Gelman.Carlin.95},
but the Jeffreys' prior has some minimax properties
\cite{Barron94} that make it more robust.
\end{itemize}
The hidden or latent variables here are the component scores
$\vec{l}$.  The model parameters are the gamma parameters
$\vec\beta$ and $\vec\alpha$, and the component loading matrix
${\mat\Theta}$. Denote the model as GP standing for Gamma-Poisson.
The full likelihood for each document, $p(\vec{w}, \vec{l} \,|\,
\vec\beta, \vec\alpha, \mat\Theta, K, \mbox{GP})$, is composed of
two parts. The first part comes from $K$ independent gamma
distributions for $l_k$, and the second part comes from $J$
independent Poisson distributions with parameters $\sum_k
l_k\theta_{j,k}$. \begin{widetext}
\begin{equation}
\label{eq-lf}
\begin{array}{cc}\textrm{likelihood of $\vec{l}$} &
\textrm{likelihood of $\vec{w}$ given $\vec{l}$}\\
\overbrace{\prod_k
\frac{\beta_k^{\alpha_k} l_k^{\alpha_k-1}
\exp\{- \beta_k l_k\}}{\Gamma(\alpha_k)}} & \overbrace{\prod_{j} \frac{
\left(\sum_k l_k\theta_{j,k}\right)^{w_j}\exp\left\{-\left(\sum_k
l_k\theta_{j,k}\right)\right\}} {w_{j}! }}
\end{array}
\end{equation}
\end{widetext}

\subsection{The Conditional Gamma-Poisson Model}

\comaOff{Why conditional GP model anyway?? example? justification?}
In practice, when fitting the parameters $\vec{\alpha}$
in the GP or DM model, it is often the case that the $\alpha_k$
go very small.
Thus, in this situation, perhaps 90\% of the component scores
$l_k$ are negligible, say less than $10^{-8}$ once normalised.
Rather than maintaining these negligible values,
we can allow component scores to be zero with some finite probability.
The Conditional Gamma-Poisson Model, denoted
CGP for short, introduces this capability.
In retrospect, CGP is a
sparse GP with an additional parameter per component to encourage
sparsity.

The CGP model extends the Gamma-Poisson model by making the $l_k$ zero
sometimes. In the general case, the $l_k$ are independent and zero
with probability $\rho_k$ and otherwise gamma distributed with
probability $1-\rho_k$.
\[
l_k ~\sim~\mbox{Gamma}(\alpha_k,\beta_k) ~~~~~~~~~~~\mbox{for
$k=1,\ldots,K$}.
\]
Denote the model as CGP standing for Conditional Gamma-Poisson,
and the full likelihood is now $p(\vec{w}, \vec{l} \,|\,
\vec\beta, \vec\alpha, \vec\rho, {\mat\Theta}, K, \mbox{CGP})$. The
 full likelihood for each document, modifying
the above Equation~(\ref{eq-lf}), replaces the term inside
$\prod_k$ with
\begin{equation}
\label{eq-clf} (1-\rho_k)\frac{\beta_k^{\alpha_k} l_k^{\alpha_k-1}
\exp\{- \beta_k l_k\}}
        {\Gamma(\alpha_k)}
+ \rho_k 1_{l_k=0}
\end{equation}

\subsection{The Dirichlet-Multinomial Model}
\label{sct-dd}

\comaOff{what's the practical consequence of not modelling $L$?}
\comaWOff{Good question.  Dont know.}

The Dirichlet-multinomial form of DCA was introduced as MPCA. In
this case, the normalised latent variables $\vec{m}$ are used, and
the total word count $L$ is not modelled.
\[
\vec{m} ~\sim~ \mbox{Dirichlet}_K(\vec\alpha) ~, ~~~~~~~~~ \vec{w}
~\sim~ \mbox{Multinomial}(L, {\mat\Theta}\vec{m})
\]
The first argument to the multinomial is the total count,
the second argument is the vector of probabilities.
Denote the model as DM, and the full likelihood is now $p(\vec{w},
\vec{m} \,|\, L, \vec\alpha,  {\mat\Theta},K, \mbox{DM})$. The
full likelihood for each document becomes:
\begin{equation}
\label{eq-slfm} C^L_{w_1, \ldots , w_J} \Gamma\left(\sum_k\alpha_k\right)
        \prod_k \frac{ m_k^{\alpha_k-1}}{\Gamma(\alpha_k)}
\prod_{j} \left(\sum_k m_k\theta_{j,k}\right)^{w_j}
\end{equation}
where $C^L_{\vec{w}}$ is $L$ choose $w_1, \ldots , w_J$. This model can also
be derived from the Gamma-Poisson model, shown in the next section.

\subsection{A Multivariate Version}
\label{sct-mv}

Another variation of the methods is to allow grouping of the count
data.  Words can be grouped into separate variable sets.  These
groups might be ``title words,'' ``body words,'' and ``topics'' in
web page analysis or ``nouns,'' ``verbs'' and ``adjectives'' in
text analysis. The groups can be treated with separate discrete
distributions, as below. The $J$ possible word types in a document
are partitioned into $G$ groups $B_1,\ldots,B_G$.  The total word
counts for each group $g$ is denoted $L_g = \sum_{j\in B_g} w_j$.
If the vector $\vec{w}$ is split up into $G$ vectors $\vec{w}_g=\{
w_j\,:\,j\in B_g\}$, and the matrix ${\mat\Theta}$ is now
normalised by group in each row, so $\sum_{j\in B_g} \theta_{j,k}
= 1$, then a multivariate version of DCA is created so that for
each group $g$,
\[
\vec{w}_g ~\sim~ \mbox{Multinomial}
        \left( L_g,\, \left\{ \sum_k m_k \theta_{j,k}\,:\,j \in B_g \right\}
    \right) ~.
\]
Fitting and modelling methods for this variation are related to
LDA or MPCA, and will not be considered in more detail here. This
has the advantage that different kinds of words have their own
multinomial and the distribution of different kinds is ignored.
This version is demonstrated subsequently on US Senate voting
records, where each multinomial is now a single vote for a particular senator.

\section{Related Work}
\label{sct-hn}

These sections begins by relating the main approaches
to each other, then placing them
in the context of exponential family models,
and finally a brief history is recounted.

\subsection{Correspondences}

Various published cases of DCA can be represented in terms of this
format, as given in Table~\ref{tbl-exp}.
A multinomial with total count $L$ and $J$ possible outcomes
is the
bagged version of $L$ discrete distributions with $J$ possible outcomes.
In the table, {\it NA}
indicates that this aspect of the model was not required
to be specified because the methodology made no use of it.
\begin{table}
\caption{Previously Published Models}
\label{tbl-exp}
\begin{center}
\begin{tabular}{|c|c|c|c|c|}
 \hline
Name & Bagged & Components & $p(\vec{x}/\vec{w}\,|\,{\mat\Theta},\vec{l})$
         & $p(\vec{l}/\vec{m})$\\
\hline
NMF \cite{lee99learning} & yes & $\vec{l}$ & Poisson &  {\it NA}\\
PLSI \cite{hofmann99probabilistic} &
      no & $\vec{m}$ & discrete &  {\it NA} \\
LDA \cite{Blei03} & no & $\vec{m}$ & discrete & Dirichlet  \\
MPCA \cite{Buntine.mpca} & yes & $\vec{m}$ & multinomial & Dirichlet  \\
GaP \cite{canny04} & yes & $\vec{l}$ & Poisson & gamma \\
\hline
\end{tabular}
\end{center}
\end{table}
\comaOff{where are the $PD_*$ things defined? why not put them in?}
Note that NMF used a cost function formulation, and thus
avoided defining likelihood models.
It is shown later that its cost function corresponds to
a Gamma-Poisson with parameters
$\vec\alpha=\vec\beta=\vec{0}$ (i.e., all zero).

LDA has the multinomial of MPCA replaced by a sequence of discrete
distributions, and thus the choose term drops, as per
Section~\ref{sct-bs}. PLSI is related to LDA but lacks a prior
distribution on $\vec{m}$.  It does not model these latent variables
using full probability theory, but instead using a weighted
likelihood method \cite{hofmann99probabilistic}. Thus PLSI is a
non-Bayesian version of LDA, although its weighted likelihood method
means it accounts for over-fitting in a principled manner.

LDA and MPCA also have a close relationship to GaP (called GP here).
If the parameter $\vec\alpha$ is treated as known
and not estimated from the data, and
the $\vec \beta$ parameter vector has the same value for each $\beta_k$,
 then $L$ is {\it aposteriori}
independent of  $\vec{m}$ and $\mat\Theta$. In this context
LDA, MPCA and GaP are equivalent models ignoring representational issues.
\begin{lemma}
Given a Gamma-Poisson model of Section~\ref{sct-gp} where the
$\vec \beta$ parameter is a constant vector with all entries
the same, $\beta$, the model is equivalent to a
Dirichlet-multinomial model of Section~\ref{sct-dd} where
$m_k=l_k/\sum_k l_k$, and in addition
\[
L ~\sim~ \mbox{Poisson-Gamma}
       \left(\sum_k \alpha_k,\beta,1\right)
\]
\end{lemma}
\begin{proof}
{\small
Consider the Gamma-Poisson model. The sum
$L=\sum_j w_j$ of Poisson variables $\vec{w}$ has the
distribution of a Poisson with parameter given by the sum of their means.
When the sum of Poisson variables is known, the set of Poisson variables
has a multinomial distribution conditioned on the sum (the total count)
\cite{Ross.intro}. The Poisson distributions on $\vec{w}$ then
is equivalent to:
\[
\begin{aligned}
L &~\sim~ \mbox{Poisson}
       \left(\sum_k l_k\right)\\
\vec{w} &~\sim~ \mbox{Multinomial}
        \left(L,\, \frac{1}{\sum_k l_k  } \mat\Theta \vec{l}
        \right)~.
\end{aligned}
\]
Moreover, if the $\vec \beta$ parameter is constant, then
$m_k=l_k/\sum_k l_k$ is distributed as
$\mbox{Dirichlet}_K (\vec\alpha)$,
and $\sum_k l_k$ is distributed independently as a
$\mbox{Gamma}(\sum_k \alpha_k,\beta)$.
The second distribution above can then be
represented as
\[
\vec{w} ~\sim~ \mbox{Multinomial}
        \left(L,\, \mat\Theta \vec{m}\right) ~.
\]
Note also, that marginalising out $\sum_k l_k$ convolves a
Poisson and a gamma distribution to produce a Poisson-Gamma
distribution for $L$ \cite{Bernardo.Smith.94}.
}
\end{proof}
If
$\vec\alpha$ is estimated from the data in GaP, then the presence
of the observed $L$ will influence $\vec\alpha$, and thus the other
estimates such as of $\mat\Theta$. In this case, LDA and MPCA
will no longer be effectively equivalent to GaP. Note, Canny
recommends fixing $\vec\alpha$ and estimating $\vec \beta$ from
the data \cite{canny04}.

To complete the set of correspondences, note that in
Section~\ref{sct-mf} it is proven that NMF corresponds to a
maximum likelihood version of GaP, and thus it also corresponds to
a maximum likelihood version of LDA, MPCA, and PLSI.

\subsection{Notes on the Exponential Family}

For the general DCA model
of Section~\ref{sct-gdca},
when $p(\vec{w}\,|\, {\mat\Theta}\vec{l} )$ is in the so-called exponential
family distributions \cite{Gelman.Carlin.95}, the expected value
of $\vec{w}$ is referred to as the {\it dual parameter}, and it is
usually the parameter we know best. For the Bernoulli with
probability $p$, the dual parameter is $p$, for the Poisson with
rate $\lambda$, the dual parameter is $\lambda$, and for the
Gaussian with mean $\mu$, the dual parameter is the mean. Our
formulation, then, can be also be interpreted as letting $\vec{w}$
be exponential family with dual parameter given by $({\mat\Theta}
\vec{l})$. Our formulation then generalises PCA in the same way
that a linear model \cite{mccullagh.nelder:glm} generalises linear
regression.

Note, an alternative has also been presented
\cite{sanjoy-generalization} where $\vec{w}$ has an exponential
family distribution with natural parameters given by $(\mat\Theta
\vec{l})$. For the Bernoulli with probability $p$, the natural
parameter is $\log (p/(1-p))$, for the Poisson with rate $\lambda$,
the natural parameter is $\log \lambda$ and for the Gaussian with
mean $\mu$, the natural parameter is the mean. This formulation
generalises PCA in the same way that a generalised linear model
\cite{mccullagh.nelder:glm} generalises linear regression.

\subsection{Historical notes}

Several independent groups within the statistical computing and
machine learning community have contributed to the development of
the DCA family of methods. Some original research includes the
following: grade of membership (GOM) \cite{GOM}, probabilistic
latent semantic indexing (PLSI) \cite{hofmann99probabilistic},
non-negative matrix factorisation (NMF) \cite{lee99learning},
genotype inference using admixtures \cite{pritchard:1999}, latent
Dirichlet allocation (LDA) \cite{Blei03}, and Gamma-Poisson models
(GaP) \cite{canny04}. Modifications and algorithms have also been
explored as multinomial PCA (MPCA) \cite{Buntine.mpca} and multiple
aspect modelling \cite{minka-uai02}. Several of these models have
been interpreted as instances of a more general family of
mixed-membership models \cite{Erosheva04}.

The first clear enunciation of the large-scale model in its
Poisson form comes from \cite{lee99learning}, and in its
multinomial form from \cite{hofmann99probabilistic} and
\cite{pritchard:1999}. The first clear expression of the problem
as a latent variable problem is given by  \cite{pritchard:1999}.  The
relationship between LDA and PLSI and that NMF was a Poisson
version of LDA was first pointed out by \cite{Buntine.mpca},
and proven in \cite{Gaussier2005}. The
connections to ICA come from \cite{BuntineJakulin04} and
\cite{canny04}. The general Gamma-Poisson formulation, perhaps the
final generalisation to this line of work, is in \cite{canny04}.

Related techniques in the statistical community can be
traced back to Latent Class Analysis 
developed in the 1950's, and a rich theory has since developed
relating the methods to correspondence analysis and other
statistical techniques \cite{HeijdenGilula99}.

\section{Component Assignments for Words}
\label{sct-caw}

In standard mixture models, each document in a collection is
assigned to one latent component. The DCA family of models can be
interpreted as making each word in each document be
assigned to one latent component.
To see this, we introduce another latent vector which represents
the component assignments for different words. As in
Section~\ref{sct-bs}, this can be done using a bag of components
or a sequence of components representation, and no effective
change occurs in the basic models, or in the algorithms so
derived. What this does is expand out the term
$\mat\Theta \vec{l}$ into parts, treating it as if it is the
result of marginalising out some latent variable.

We introduce a $K$-dimensional discrete latent vector $\vec{c}$
whose total count is $L$, the same as the word count. The
count $c_k$ gives the number of words in the document appearing in
the $k$-th component.  Its posterior mean makes a good
diagnostic and interpretable result. A document from the sports
news might have 50 ``football'' words, 10 ``German'' words, 20
``sports fan'' words and 20 ``general vocabulary'' words.
\comaOff{It's perhaps important that you explain the relation between
V and L better. And I can't stress how confusing it is to mix
vectors and matrices. I.e., I'm still not sure why matrix
$\vec{V}$ would be different from $\vec{L}$.}
\comaOff{OK, general mess with the notation. It would be a good thing
to use normal matrix notation throughout the paper, it's messy as
it is now. So, a vector is $\vec{x}$, a matrix is $\vec{X}$, a
column from that matrix is $(\vec{X})_{j}$, and a row from that
matrix is $\vec{x}_{j}$? }
\comaWOff{OK, largely adopted except for the column vector.}

This latent vector is derived from a larger latent matrix,
$\mat{V}$  of size $J\times K$ and entries $v_{j,k}$.
This has row totals $w_j$ as given in the observed
data and column totals $c_k$.  Vectors $\vec{w}$ and $\vec{c}$ are these
word appearance counts and component appearance counts, respectively,
based on summing rows and columns of matrix $\mat{V}$.
This is shown in Figure~\ref{fig-V}.
\begin{figure}
\begin{center}
\[
\begin{array}{ccc}
 & {K \textrm{ components}} \\
 \rotatebox{90}{\hspace{-5mm}J \textrm{ words}} &
\hspace{6mm} \underline{\overbrace{\hspace{-6mm}\left\{
\left(\begin{array}{cccc}
v_{1,1} & v_{1,2} & \cdots & v_{1,K}\\
v_{2,1} & v_{2,2} & \cdots & v_{2,K}\\
\\
\vdots & \vdots & \ddots & \vdots\\
\\
v_{J,1} & v_{J,2} & \cdots & v_{J,K}\\
\end{array}\right)\right|\hspace{-4mm}}} \hspace{4mm}
&\hspace{-1mm}
\begin{array}{c}
w_1\\
w_2\\
\\
\vdots\\
\\
w_J
\end{array}\\
& ~~~\begin{array}{cccc} c_{1}~~ & c_{2}~~ & \cdots &
c_{K}~~\end{array} &
\end{array}
\]
\end{center}
\caption{A representation of a document as a
contingency table.}
\label{fig-V}
\end{figure}

The introduction of the latent matrix $\mat{V}$ changes the forms of
the likelihoods, and makes the development and analysis of
algorithms easier. This section catalogues the likelihood formula,
to be used when discussing algorithms. The choices of different
statistical word count models have been recently evaluated by
\cite{airoldi2005bayesian} in the context of document
classification.

\comaWOff{Not sure this next stuff is used}

\subsection{The Gamma-Poisson Model}

With the new latent matrix $\mat{V}$, the distributions underlying the
Gamma-Poisson model become:
\begin{eqnarray}
\label{eq-gpc} l_k &~\sim~&
   \mbox{Gamma}(\alpha_k,\beta_k) \\
\nonumber
c_{k} &~\sim~& \mbox{Poisson}\left( l_k \right)~\\
\nonumber w_j& = &\sum_k v_{j,k}, ~~\mbox{where}~ v_{j,k} \sim
\mbox{Multinomial}\left(c_k, \vec\theta_{\cdot,k}\right)~.
\end{eqnarray}

The joint likelihood for a document, $p({\mat{V}}, \vec{l} \,|\, \vec\beta,
\vec\alpha, \mat \Theta, K, \mbox{GP})$ (the $\vec{w}$ are now
derived quantities so not represented), thus becomes, after some
rearrangement
\begin{equation}
\label{eq-gpl1} \prod_k \frac{\beta_k^{\alpha_k}
l_k^{c_k+\alpha_k-1} \exp\{- (\beta_k+1) l_k\}}
         {\Gamma(\alpha_k)}
\prod_{j,k} \frac{ \theta_{j,k}^{v_{j,k}} } {v_{j,k}! }~.
\end{equation}
Note that $\vec{l}$ can be marginalised out, yielding
\begin{equation}
\label{eq-gpl2} \prod_k
\frac{\Gamma(c_k+\alpha_k)}{\Gamma(\alpha_k)}
        \frac{\beta_k^{\alpha_k}}{(\beta_k+1)^{c_k+\alpha_k}}
\prod_{j,k} \frac{ \theta_{j,k}^{v_{j,k}} } {v_{j,k}! }~.
\end{equation}
and the posterior mean of $l_k$ given $\vec{c}$ is
$(c_k+\alpha_k)/(1+\beta_k)$. Thus each $c_k \sim
\mbox{Poisson-Gamma}(\alpha_k,\beta_k,1)$.

\subsection{The Conditional Gamma-Poisson Model}

The likelihood follows the GP case,
except that with probability $\rho_k$, $l_k=0$ and
thus $c_k=0$.
The joint likelihood, $p(\mat{V}, \vec{l} \,|\, \vec\beta,
\vec\alpha, \vec\rho, \mat\Theta,K, \mbox{CGP})$,
 thus becomes, after some rearrangement
\begin{eqnarray}
\label{eq-gpcl1}\begin{split} \prod_k \left( (1-\rho_k)\left(
\frac{\beta_k^{\alpha_k} l_k^{c_k+\alpha_k-1} \exp\{- (\beta_k+1)
l_k\}}
         {\Gamma(\alpha_k)}
\prod_{j} \frac{ \theta_{j,k}^{v_{j,k}} } {v_{j,k}! } \right)\right. \\
\nonumber  \left. + \rho_k \left( 1_{l_k=0}  1_{c_k=0} \prod_j
1_{v_{j,k}=0} \right) \right) ~.
\end{split}
\end{eqnarray}
Note that $\vec{l}$ can be marginalised out, yielding
\begin{equation}
\begin{split} \label{eq-gpcl2} \prod_k \left(
 (1-\rho_k)
\frac{\Gamma(c_k+\alpha_k)}{\Gamma(\alpha_k)}
        \frac{\beta_k^{\alpha_k}}{(\beta_k+1)^{c_k+\alpha_k}}
+ \rho_k   1_{c_k=0} \right) \\\prod_{j} \frac{
\theta_{j,k}^{v_{j,k}} } {v_{j,k}! } ~.\end{split}
\end{equation}
The $\theta_{j,k}$ can be pulled out under the constraint $\sum_j
v_{j,k}=c_k$. The posterior mean of $l_k$ given $\vec{c}$ is
$(1-\rho_k)(c_k+\alpha_k)/(1+\beta_k)$.

\subsection{The Dirichlet-Multinomial Model}

For the Dirichlet-multinomial model, a similar reconstruction
applies:
\begin{eqnarray}
\label{eq-ddc}
\vec{m} &~\sim~& \mbox{Dirichlet}_K(\vec\alpha) ~\\
\nonumber
c_{k} &~\sim~& \mbox{Multinomial}\left(L, \vec{m} \right)~\\
\nonumber w_j& = &\sum_k v_{j,k},~ \mbox{where}~ v_{j,k} \sim
\mbox{Multinomial}\left(c_k, \vec\theta_{\cdot,k}\right)~.
\end{eqnarray}
The joint likelihood, $p({\mat{V}}, \vec{m} \,|\, \vec\alpha, \mat\Theta, K,
\mbox{DM})$,
 thus becomes, after some rearrangement
\begin{equation}
\label{eq-ddl1} L! \,\Gamma\left(\sum_k \alpha_k\right) \prod_k
\frac{m_k^{c_k+\alpha_k-1}} {\Gamma(\alpha_k)} \prod_{j,k} \frac{
\theta_{j,k}^{v_{j,k}}} {v_{j,k}! } ~.
\end{equation}
Again, $\vec{m}$ can be marginalised out yielding
\begin{equation}
\label{eq-ddl2} L! \,\frac{\Gamma\left(\sum_k \alpha_k\right)}
      {\Gamma\left(L+\sum_k \alpha_k\right)}
\prod_k \frac{\Gamma(c_k+\alpha_k)} {\Gamma(\alpha_k)} \prod_{j,k}
\frac{ \theta_{j,k}^{v_{j,k}}} {v_{j,k}! } ~.
\end{equation}

\section{Algorithms}
\label{sct-alg}

In developing an algorithm, the standard approach is to
match an optimization algorithm to the functional form of
the likelihood.  When using Bayesian or some other
statistical methodology, this basic approach is
usually a first step, or perhaps an inner loop for some
more sophisticated computational statistics.

The likelihoods do not yield easily to standard EM analysis.
To see this, consider the forms of the likelihood for a single
document for the GP model, and consider the probability
for a latent variable $\vec{z}$ given the observed data $\vec{w}$,
$p(\vec{z} \,|\, \vec{w}, \vec\beta, \vec\alpha, \mat\Theta,K,
\mbox{GP})$.
For EM analysis, one needs to be able to compute
$\mathbb{E}_{\vec{z} \sim p(\vec{z}|\vec{w},\mat \Theta, \ldots )}\left[
 \log p(\vec{w},\vec{z} \,|\,  \mat \Theta, \ldots )\right]$.
There are three different forms of the likelihood seen so far
depending on which latent variables $\vec{z}$ are kept on
the left-hand side of the probability:
\begin{description}
\item[$p(\vec{w}, \vec{l} \,|\, \vec\beta, \vec\alpha, \mat\Theta,K,
\mbox{GP})$: ]
from Equation~(\ref{eq-lf}) has the term
$\left(\sum_k l_k\theta_{j,k}\right)^{w_j}$,
which means there is no known simple posterior distribution
for $\vec{l}$ given $\vec{w}$.
\item[$p(\vec{w}, \vec{l}, \mat V \,|\, \vec\beta, \vec\alpha, \mat\Theta,K,
\mbox{GP})$: ]
from Equation~(\ref{eq-gpl1}) has the term $l_k^{c_k+\alpha_k-1}$
which links the two latent variables $\vec{l}$
and $\mat V$,
and prevents a simple evaluation of
$\mathbb{E}_{\vec{l}, \mat V}[ v_{j,k}]$ as required
for the expected log probability.
\item[$p(\vec{w}, \mat V \,|\, \vec\beta, \vec\alpha, \mat\Theta,K,
\mbox{GP})$: ]
from Equation~(\ref{eq-gpl2})
has the term $\Gamma(c_k+\alpha_k)$ (where $c_k=\sum_j v_j,k$), which
means there is no known simple posterior distribution for $\mat V$
given $\vec{w}$.
\end{description}
\noindent Now one could always produce an EM-like algorithm by
separately updating $\vec{l}$ and $\mat V$ in turn according to some
mean formula, but the guarantee of convergence of $\mat \Theta$ to a
maximum posterior or likelihood value will not apply. In this spirit
earlier authors point out that EM-like principles apply and use EM
terminology since EM methods would apply if $\vec{l}$ was observed
\footnote{The likelihood $p(\vec{w}, \mat V \,|\, \vec{l},
\vec\beta, \vec\alpha,\mat\Theta, K, \mbox{GP})$ can be treated with
EM methods using the latent variable $ \mat V$ and leaving $\vec{l}$
as if it was observed.}. For the exponential family, which this
problem is in, the variational approximation algorithm with
Kullback-Leibler divergence corresponds to an extension of the EM
algorithm \cite{ghahramani00propagation, Buntine.mpca}. This
variational approach is covered below.

Algorithms for this problem follow some general approaches in the
statistical computing community.
Three basic approaches are presented here: a
 variational approximation,
Gibbs sampling, and Rao-Blackwellised Gibbs sampling.
A maximum likelihood algorithm is not presented
because it can be viewed as a simplification of the algorithms here.

\subsection{Variational Approximation with Kullback-Leibler Divergence} \label{sct-mf}

This approximate method was first applied to the sequential variant
of the Dirichlet-multinomial version of the problem by
\cite{Blei03}.
A fuller treatment of these variational methods for the
exponential family is given in
\cite{ghahramani00propagation,Buntine.mpca}.

In this approach a factored
posterior approximation is made for the latent variables:
\[
p(\vec{l}, \mat V \,|\, \vec{w}, \vec\beta, \vec\alpha, \mat\Theta,K,
\mbox{GP})
~\approx ~q(\vec{l},\mat V)~= ~
q_{\vec{l}} (\vec{l})q_{\mat V} (\mat V)
\]
and this approximation is used to find expectations as part of an
optimization step.
The EM algorithm results if an equality holds.
The functional form of the
approximation can be derived by inspection of the recursive
functional forms (see \cite{Buntine.mpca}
Equation~(\ref{eq-slfm})):
\begin{eqnarray}
\label{eq-ff} \nonumber q_{\vec{l}}(\vec{l}) & \propto & \exp \left(
\mathbb{E}_{\mat V \sim q_{\mat V}({\mat{V}})}\left[
   \log p\left(
      \vec{l},{\mat{V}},\vec{w} \,| \, {\mat\Theta},\vec{\alpha},\vec{\beta},K
   \right) \right]
\right)\\
q_{\mat V}({\mat{V}}) & \propto & \exp \left(
\mathbb{E}_{\vec{l}\sim q_{\vec{l}}(\vec{l})} \left[ \log
p\left(\vec{l},{\bf v},\vec{w} \,| \, {\mat
\Theta},\vec{\alpha},\vec{\beta},K\right) \right] \right)~.
\end{eqnarray}

An important computation used during convergence in this approach
is a lower bound on the individual document log probabilities.
This naturally
falls out during computation (see \cite{Buntine.mpca}
Equation~(\ref{eq-gpl1})).
Using the approximation $q(\vec{l},\mat{V})$ defined
by the above proportions, the bound is given by
\begin{eqnarray*}
\begin{split}\log p\left(\vec{w} \,| \, {\bf
\Theta},\vec{\alpha},\vec{\beta},K\right) \geq \\ \label{eq-pa}
\mathbb{E}_{\vec{l},{\mat{V}} \sim q(\vec{l},\mat{V})}
  \left[ \log \,
     p\left(\vec{l},{\mat{V}},
     \vec{w}
      \,| \, {\mat \Theta},\vec{\alpha},\vec{\beta},K\right)
      \right]\\
 + I(q_{\vec{l}}(\vec{l})) + I(q_{\mat{V}}(\mat{V}))~.
 \end{split}
\end{eqnarray*}
The variational approximation applies to the Gamma-Poisson version and
the Dirichlet-multinomial version.

\comaOff{OK, the bound in freaky notation is just
saying that KL-divergence $D(P(\vec{l},\vec{v}|{\cal D}ata)
 \|q(\vec{l}|{\cal D}ata)q(\vec{v}|{\cal D}ata))$ bounds the
decrease in likelihood per a document or something?}

\subsubsection{For the Gamma-Poisson Model:}
Looking at the recursive functionals of Equation\ (\ref{eq-ff})
and the likelihood of Equation~(\ref{eq-gpl1}),
it follows that
$q_{\vec{l}}()$ must be $K$ independent Gammas one for each component, and
$q_{\mat V}()$ must be $J$ independent multinomials,
one for each word.
The most general case for the approximation $q()$ is thus
\begin{eqnarray*}
l_k & \sim & \mbox{Gamma}(a_k,b_k) \\
\{v_{j,1\ldots K}\} & \sim &
\mbox{Multinomial}(w_j,\{n_{j,k}:k=1,\ldots,K\})~,
\end{eqnarray*}
which uses approximation parameters
$(a_k,b_k)$ for each Gamma and
and $\vec{n}_{\cdot,k}$ (normalised as $\sum_k n_{j,k} = 1$)
for each multinomial.
These parameters form two vectors $\vec{a},\vec{b}$
and a matrix $\mat{N}$ respectively.
The approximate posterior takes the form
$q_{\vec{l}}(\vec{l}\,|\,\vec{a},\vec{b})q_{\mat V} (\mat V\,|\,\mat{N})$.

Using these approximating distributions, and
again looking at the recursive functionals of Equation\ (\ref{eq-ff}),
one can extract the rewrite rules for the parameters:
\begin{eqnarray}
\label{eq-grw}
n_{j,k} &= & \frac{1}{Z_j} \theta_{j,k} \exp\left( \mathbb{E}\left[\log\, l_k \right] \right) ~,\\
\nonumber a_k &=& \alpha_k + \sum_j w_j n_{j,k}  ~,\\
\nonumber b_k &=& 1 + \beta_k  ~,\\
\nonumber \mbox{\bf where } \mathbb{E}\left[\log\, l_k \right]
&\equiv& \mathbb{E}_{l_k \sim p(l_k\,|\, a_k.b_k)} \left[\log\, l_k
\right] \\ && ~=~
\Psi_0(a_k) - \log b_k ~,\\
\nonumber Z_j &\equiv& \sum_k \theta_{j,k}\exp\left( \mathbb{E}\left[\log\, l_k \right] \right) ~.
\end{eqnarray}
Here, $\Psi_0()$ is the digamma function, defined
as $\frac{\textrm{d} \ln \Gamma(x)}{\textrm{d} x}$
and available in most scientific libraries.
These equations form the first step of each major cycle, and are
performed on each document.

The second step is to re-estimate the model parameters
$\mat\Theta$ using the posterior approximation by maximising the
expectation of the log of the full posterior probability
\[
\mathbb{E}_{\vec{l},{\mat{V}} \sim q_{\vec{l}} (\vec{l})q_{\mat V} (\mat V)}
       \left[\log \,
     p\left(\vec{l},{\mat{V}},
     \vec{w},{\mat\Theta}
      \,| \, \vec{\alpha},\vec{\beta},K\right) \right] ~.
\]
This incorporates Equation~(\ref{eq-gpl1}) for each document, and
a  prior for each $k$-th column of
${\mat\Theta}$ of $\mbox{Dirichlet}_J(\vec\gamma)$
(the last model item in Section~\ref{sct-gp}). Denote the intermediate
variables $n_{j,k}$ for the $i$-th document by adding a $(i)$
subscript, as $n_{j,k,(i)}$, and likewise for $w_{j,(i)}$. All
these log probability formulas yield linear terms in
$\theta_{j,k}$, thus with the normalising constraints for
$\mat\Theta $ one gets
\begin{equation}
\label{eq-grw2} \theta_{j,k} ~\propto  ~ \sum_i w_{j,(i)}
n_{j,k,(i)} + \gamma_j ~.
\end{equation}
\noindent The lower bound on the log probability of
Equation~(\ref{eq-pa}), after some simplification and use of the
rewrites of Equation~(\ref{eq-grw}), becomes
\begin{equation}
\begin{split}
\label{eq-grw3} \log\frac{1}{\prod_j w_j!}
 -  \sum_k \log \frac{\Gamma(\alpha_k)b_k^{a_k}}{\Gamma(a_k)
 \beta_k^{\alpha_k}}
 + \sum_k (\alpha_k-a_k)\mathbb{E}\left[\log\, l_k \right]\\
+ \sum_j w_j \log Z_j ~.
\end{split}
\end{equation}
The variational approximation algorithm
for the Gamma-Poisson version is summarised in
Figure~\ref{fig-mfa}.
An equivalent algorithm is produced if words are presented sequentially
instead of being bagged.

\begin{figure*}[htb]
\begin{center}
\begin{boxedminipage}{6.8in}
\begin{enumerate}
\item Initialise $\vec{a}$ for each document. The uniform
initialisation would be $a_k=\left(\sum_k \alpha_k + L\right)/K$.
Note $\mat{N}$ is not stored.

\item Do for each document:
\begin{enumerate}
\item Using Equations~(\ref{eq-grw}),
recompute $\mat{N}$ and update $\vec{a}$ in place.

\item Concurrently, compute the log-probability bound of
Equation~(\ref{eq-grw3}),
and add to a running total. \item Concurrently, maintain the
sufficient statistics for ${\mat\Theta}$, the total $\sum_i
w_{j,(i)} n_{j,k,(i)}$ for each $j,k$ over documents. \item Store
$\vec{a}$ for the next cycle and discard $\mat{N}$.
\end{enumerate}

\item Update ${\mat\Theta}$ using Equation~(\ref{eq-grw2}),
normalising appropriately. \item Report the total log-probability
bound, and repeat, starting at Step~2.
\end{enumerate}
\end{boxedminipage}
\end{center}

\caption{K-L Variational Algorithm for Gamma-Poisson} \label{fig-mfa}
\end{figure*}

\paragraph{Complexity: }
Because Step~2(a) only uses words appearing in a document, the
full Step~2 is $O(SK)$ in time complexity where $S$ is the number
of words in the full collection. Step~3 is $O(JK)$ in time
complexity. Space complexity is $O(IK)$ to store the intermediate
parameters $\vec{a}$ for each document, and the $O(2JK)$ to store
${\mat\Theta}$ and its statistics. In implementation, Step~2 for
each document is often quite slow, and thus both  $\vec{a}$ and
the document word data can be stored on disk and streamed, thus
the main memory complexity is $O(2JK)$ since the
$O(S)$ and $O(IK)$ terms are on disk.
If documents are very small (e.g., $S/I \ll K$, for instance
``documents'' are sentences or phrases),
then this does not apply.

\paragraph{Correspondence with NMF:}
\label{sct-nmf}
A precursor to the GaP model is non-negative matrix factorisation
(NMF) \cite{lee99learning}, which is based on the matrix
approximation paradigm using Kullback-Leibler divergence. The
algorithm itself, converted to the notation used here, is as
follows

\begin{align*} l_{k,(i)} \, &\longleftarrow \,
      l_{k,(i)} \sum_j \frac{\theta_{j,k}}{\sum_j \theta_{j,k}}
                       \frac{ w_{j,(i)} }
                        { \sum_k \theta_{j,k} l_{k,(i)} }
\\ \theta_{j,k} \, &\longleftarrow \,
      \theta_{j,k} \sum_i \frac{l_{k,(i)}}{\sum_i l_{k,(i)}}
                       \frac{ w_{j,(i)} }
                        { \sum_k \theta_{j,k} l_{k,(i)} }\end{align*}

Notice that the solution is indeterminate up to a factor $\psi_k$.
Multiply $l_{k,(i)}$ by $\psi_k$ and divide $\theta_{j,k}$ by
$\psi_k$  and the solution still holds. Thus, without loss of
generality, let $\theta_{j,k}$ be normalised on $j$, so that $\sum_j
\theta_{j,k}=1$.
\begin{lemma}
The NMF equations above,
where ${\mat\Theta}$ is returned normalised, occur at a maxima w.r.t.\
${\mat\Theta}$ and $\vec{l}$ for the Gamma-Poisson likelihood
$\prod_i p(\vec{w}_{(i)}\,|\, {\mat\Theta}, \vec{l}_{(i)},\vec\alpha=0, \vec\beta=0,
K,\mbox{GP})$.
\end{lemma}
\begin{proof}
{\small To see this, the following will be proven. Take a solution
to the NMF equations, and divide $\theta_{j,k}$ by a factor
$\psi_k=\sum_j \theta_{j,k}$, and multiply $l_{k,(i)}$ by the same
factor. This is equivalent to a solution for the following rewrite
rules
\begin{align*} l_{k,(i)} \, &\longleftarrow \,
      l_{k,(i)} \sum_j \theta_{j,k} \frac{ w_{j,(i)} }
                        { \sum_k \theta_{j,k} l_{k,(i)} }
\\\theta_{j,k} \, &\propto \,
      \theta_{j,k} \sum_i l_{k,(i)}
                       \frac{ w_{j,(i)} }
                        { \sum_k \theta_{j,k} l_{k,(i)} }
\end{align*}
where $\theta_{j,k}$ is kept normalised on $j$.
These equations hold at a maxima to the likelihood
$\prod_i
p(\vec{w}_{(i)}\,|\, {\mat\Theta}, \vec{l}_{(i)},\vec\alpha=0, \vec\beta=0,
K,\mbox{GP})$.
The left equation corresponds to a maxima
w.r.t.\  $\vec{l}_{(i)}$
(note the Hessian for this is easily shown to be
negative indefinite), and the right
is the EM equations for the likelihood. w.r.t.\ ${\mat\Theta}$.

To show equivalence of the above and the NMF
equations, first prove the forward direction. Take the scaled solution to
NMF. The NMF equation for $l_{k,(i)}$ is equivalent to the
equation for $l_{k,(i)}$ in the lemma. Take the NMF equation for
$\theta_{j,k}$ and separately normalise both sides.  The $\sum_i
l_{k,(i)}$ term drops out and one is left with the equation for
$\theta_{j,k}$ in the lemma. Now prove the backward direction. It
is sufficient to show that the NMF equations hold for the solution
to the  rewrite rules in the lemma, since $\theta_{j,k}$ is
already normalised. The NMF equation for $l_{k,(i)}$ clearly
holds. Assuming the rewrite rules in the lemma hold, then
\begin{eqnarray*}
\theta_{j,k} & = &
      \frac{ \theta_{j,k} \sum_i  \left(l_{k,(i)}
                       w_{j,(i)} \left/
                        \sum_k \theta_{j,k} l_{k,(i)} \right.  \right) }
       {\sum_j \theta_{j,k} \sum_i  \left(l_{k,(i)}
                       w_{j,(i)} \left/
                        \sum_k \theta_{j,k} l_{k,(i)} \right.\right) }
\\ &=&
       \frac{ \theta_{j,k} \sum_i \left( \l_{k,(i)}
                       w_{j,(i)} \left/
                        \sum_k \theta_{j,k} l_{k,(i)}\right. \right) }
       {\sum_i  l_{k,(i)}
                    \sum_j \left( \theta_{j,k} w_{j,(i)}\left/
                        \sum_k \theta_{j,k} l_{k,(i)}\right. \right)   }
~\mbox{(reorder sum)}
\\ &=&
       \frac{ \theta_{j,k} \sum_i \left( l_{k,(i)}
                       w_{j,(i)} \left/
                        \sum_k \theta_{j,k} l_{k,(i)}\right. \right)   }
       {\sum_i  l_{k,(i)} }
~\mbox{(apply rewrite)}
\end{eqnarray*}
Thus the second equation for NMF holds.
}
\end{proof}
Note, including a latent variable such
as $\vec{l}$ in the likelihood (and not dealing with it
using EM methods) does not achieve
a correct maximum likelihood solution for
the expression
$\prod_i p(\vec{w}_{(i)}\,|\, {\mat\Theta}, \vec\alpha=0, \vec\beta=0,
K,\mbox{GP})$.
In practice, this is a common approximate method for handling
latent variable problems, and can lead more readily to over-fitting.

\comaOff{the above paragraph looks interesting, but I get lost. What
do you mean by "not correct maximum likelihood"? You mean it's not
globally maximum? What would you say about ordinary EM-estimation
of mixture models? }

\subsubsection{For the Dirichlet-Multinomial Model:}
The variational approximation  takes
a related form. The approximate posterior is given by:
\begin{eqnarray*}
\vec{m} & \sim & \mbox{Dirichlet}(\vec{a}) \\
\{v_{j,1\ldots{}K} \}& \sim & \mbox{Multinomial}(w_j,
\{n_{j,k}:k=1,\ldots,K\})
\end{eqnarray*}
This yields the same style update equations as
Equations~(\ref{eq-grw}) except that $\beta_k=1$
\begin{eqnarray}
\label{eq-prw}
n_{j,k} &= & \frac{1}{Z_j} \theta_{j,k} \exp\left( \mathbb{E}\left[\log\, m_k \right] \right)~,\\
\nonumber a_k &=& \alpha_k + \sum_j w_j n_{j,k}  ~,\\
\nonumber \mbox{\bf where } \mathbb{E}\left[\log\, m_k
\right]&\equiv& \mathbb{E}_{m_k \sim p(m_k\,|\, \vec{a})}
\left[\log\, m_k\right]\\&& =
\Psi_0(a_k)-\Psi_0\left(\sum_k a_k\right)\\
\nonumber Z_j&\equiv&   \sum_k \theta_{j,k}\exp\left(
\mathbb{E}\left[\log\, m_k \right]\right) ~.
\end{eqnarray}
Equation~(\ref{eq-grw2}) is also the same. The lower bound on the
individual document log probabilities, $\log p\left(\vec{w}\,| \,
{\mat\Theta},\vec{\alpha},K,\mbox{DM}\right)$ now takes the form
\begin{equation}
\begin{split} \label{eq-prw3} \log \left( C^L_{\vec{w}}\right) - \log
\frac{\Gamma\left(\sum_k a_k\right) \prod_k \Gamma(\alpha_k) }
        {\Gamma\left(\sum_k \alpha_k\right)\prod_k \Gamma(a_k) }\\
+ \sum_k (\alpha_k-a_k)\mathbb{E}\left[\log\, m_k \right] + \sum_j
w_j \log Z_j ~. \end{split}\end{equation}
The correspondence with
Equation~(\ref{eq-grw3}) is readily seen.

The algorithm for Dirichlet-multinomial version is
related to that in Figure~\ref{fig-mfa}. Equations~(\ref{eq-prw})
replace Equations~(\ref{eq-grw}), Equation~(\ref{eq-prw3})
replaces Equation~(\ref{eq-grw3}), and the initialisation for
$a_k$ should be 0.5, a Jeffreys prior.

\subsection{Direct Gibbs Sampling}

\comaWOff{In Gibbs sampling, the component counts $c_k$ introduced in
Section~\ref{sct-caw} are eliminated by replacing them with the
sum $\sum_j v_{j,k}$.}\comaOff{Why not do this in previous section too?}
\comaWOff{This then allows sampling to proceed under
the constraint that $w_j=\sum_k v_{j,k}$.}

There are two styles of Gibbs sampling that apply to DCA.
The first is a basic Gibbs sampling first proposed
by Pritchard, Stephens and Donnelly \cite{pritchard:1999}.
Gibbs sampling is a conceptually simple method.
Each unobserved variable in the problem is resampled in turn
according to its conditional distribution.
We compute its posterior distribution conditioned
on all other variables, and then sample a new value for the
variable using the posterior.
For instance, an ordering we might use in this problem is:
$\vec{l}_{(1)},\mat{V}_{(1)}$, $\vec{l}_{(2)},\mat{V}_{(2)}$,
\ldots,  $\vec{l}_{(I)},\mat{V}_{(I)}$, $\mat{\Theta}$.
All the low level sampling in this section
use well known distributions such as gamma or multinomial,
and are available in standard scientific libraries.

To develop this approach for the Gamma-Poisson, look at the full
posterior, which is a product of individual document likelihoods
with the prior for $\mat\Theta$ from
the last model item in Section~\ref{sct-gp}.  The constant terms
have been dropped.
\begin{equation}
\begin{split}
\label{eq-gpl-g} \prod_i \left( \prod_k \frac{\beta_k^{\alpha_k}
l_{k,(i)}^{c_{k,(i)}+\alpha_k-1} \exp\{- (\beta_k+1) l_{k,(i)}\}}
         {\Gamma(\alpha_k)} \prod_{j,k} \frac{ \theta_{j,k}^{v_{j,k,(i)}} } {v_{j,k,(i)}! }
\right)\\ \prod_{j,k} \theta_{j,k}^{\gamma_j} \end{split}
\end{equation}
Each of the conditional distributions used in the Gibbs
sampling are proportional to this.
The first conditional distribution is
$p(\vec{l}_{(i)} \,|\, \mat{V}_{(i)}, \vec\beta,
\vec\alpha, \mat \Theta, K, \mbox{GP})$.
From this, isolating the terms just in $\vec{l}_{(i)}$,
we see that each $l_{k,(i)}$ is conditionally gamma distributed.
Likewise, each
$\vec{v}_{\cdot,k,(i)}$ is multinomial distributed given
$\vec{l}_{(i)}$ and $\mat{\Theta}$, and each $\vec{\theta}_{\cdot,k}$
is Dirichlet distributed given all the
$\vec{l}_{(i)}$ and $\mat{V}_{(i)}$ for each $i$.
The other models are similar.  An additional effort is required to
arrange the parameters and sequencing for efficient use of memory.

\comaOff{OK, it's not even clear what is going on with Gibbs
sampling. Some basic explanation is needed - one paragraph or at
least a few lines. The figure is very helpful, but we should refer
the reader to it. Table 2 should be next to Fig 3. BTW, Table 2 is
a bit confusing - how does one sample from a conditional
Gamma-Poisson? Which variable is sampled?}

\comaOff{Perhaps the basic notion that needs to be explained is
`sampling'. By sampling we merely toss random numbers from a
certain distribution -- we randomly assign values to a random
variable in accordance with a particular distribution. Where does
one find how to sample from these distributions?}

The major differentiator for Gibbs sampling is the resampling of the
latent component vector $\vec{l}$.
The sampling schemes used for each version are
given in Table~\ref{tbl-gs}.
Some care is required with the
conditional Gamma-Poisson. When $c_k=0$, the sampling for $l_k$
needs to decide whether to use the zero case or the non-zero case.
This uses Equation~(\ref{eq-gpcl2}) to make the decision, and then
resorts to Equation~(\ref{eq-gpcl1}) if it is non-zero.
\begin{table}
\begin{center}
\begin{tabular}{|l|p{3in}|}
\hline Model & ~~~~  ~~~ ~~~ ~~~ Sampling\\\hline GP &
\begin{minipage}{3in}
\vspace*{3pt} $l_k \sim \mbox{Gamma}(c_k+\alpha_k,1+\beta_k)$.
\end{minipage}\\[5pt]\hline
CGP &
\begin{minipage}{3in}
\vspace*{3pt} If $c_k=0$, then Conditional Gamma-Poisson with rate
$\frac{p_k(1+\beta_k)^{\alpha_k}}{(1-p_k)\beta_k^{\alpha_k}+p_k(1+\beta_k)^{\alpha_k}}$ and $\mbox{Gamma}(\alpha_k,1+\beta_k)$.
If $c_k\neq 0$, revert to the above Gamma-Poisson case.\\
\end{minipage}\\[5pt]\hline
DM &
\begin{minipage}{3in}
\vspace*{3pt} $\vec{m} \sim
\mbox{Dirichlet}(\{c_k+\alpha_k\,:\,k\})$.
\end{minipage}\\[5pt]\hline
\end{tabular}
\end{center}
\caption{Sampling components for direct Gibbs on a single document}
\label{tbl-gs}
\end{table}

The direct Gibbs algorithm for the general case is given in
Figure~\ref{fig-ga}.
\begin{figure*}[htb]
\begin{center}
\begin{boxedminipage}{5.8in}
\begin{enumerate}
\item For each document $i$, retrieve the last $\vec{c}_{(i)}$ from store,
then
\begin{enumerate}
\item Sample the latent component variables $\vec{l}_{(i)}$
(or its normalised counterpart $\vec{m}_{(i)}$)
as per Table~\ref{tbl-gs}.
\item For each word $j$ in the document with positive count $w_{j,(i)}$,
the component counts vector, from Equation~(\ref{eq-gpc}) and
Equation~(\ref{eq-ddc}),
\[
\{v_{j,k,(i)}\,:\,k=1,\ldots,K\} ~ \sim ~ \mbox{Multinomial}\left(w_{j,(i)},
\left\{ \frac{l_{k,(i)}\theta_{j,k} }  { \sum_k l_{k,(i)}\theta_{j,k}  } \,:\,
k \right\} \right) ~.
\]
Alternatively, if the sequence-of-components version is to be
used, the component for each word can be sampled in turn using the
corresponding Bernoulli distribution. \item Concurrently,
accumulate the log-probability $p\left( \vec{w}_{(i)}\,| \, \vec{l}_{(i)},{\bf
 \vec{\alpha},\vec{\beta},\Theta}, K,\mbox{GP}\right)$,
$p\left( \vec{w}_{(i)}\,| \,  \vec{l}_{(i)},
      \vec{\alpha},\vec{\beta},\vec{\rho},{\mat\Theta},K,\mbox{CGP}\right)$,
or $p\left( \vec{w}_{(i)}\,| \, \vec{m}_{(i)},L_{(i)},
      \vec{\alpha},\vec{\beta},{\mat\Theta},K,\mbox{DM}\right)$.
\item Concurrently, maintain the sufficient
statistics for ${\mat\Theta}$, the total $\sum_i v_{j,k,(i)}$ for
each $j,k$ over documents.
\item Store $\vec{c}_{(i)}$ for the next
cycle and discard ${\mat{V}}_{(i)}$.
\end{enumerate}
\item Using a Dirichlet prior for rows of ${\mat\Theta}$, and
having accumulated all the counts ${\mat{V}}_{(i)}$ for each document in
sufficient statistics for ${\mat\Theta}$, then its posterior has
rows that are Dirichlet. Sample. \item Report the total
log-probability, and report.
\end{enumerate}
\end{boxedminipage}
\end{center}
\caption{One Major Cycle of Gibbs Algorithm for DCA}
\label{fig-ga}
\end{figure*}
This Gibbs scheme turns out to correspond to the variational approximation, excepting that sampling is done instead of
maximisation or expectation. \comaOff{this is good to say, but it's
not clear why it would be similar}

The log probability of the words $\vec{w}$ can also accumulated
in step~1(c). While they are in terms of the latent variables,
they still represent a reasonably unbiased estimate of the
likelihoods such as
$p\left( \vec{w}_{(1)},\ldots,\vec{w}_{(I)} \,| \,
  \vec{\alpha},\vec{\beta},{\mat\Theta},K,\mbox{GP}\right)$.

\subsection{Rao-Blackwellised Gibbs Sampling}

Rao-Blackwellisation of Gibbs sampling \cite{casella96}
combines closed form updates of variables with Gibbs sampling.
It does so by a process called marginalisation
or variable elimination.  When feasible, it can lead to significant
improvements, the general case for DCA.  Griffiths and Steyvers
\cite{griffiths04pnas}
introduced this algorithm for LDA, and it easily extends to
the Gamma-Poisson model and its conditional
variant with little change to the sampling routines.

When using this approach, the first step is to consider the
full posterior probability and see which variables can be
marginalised out without introducing computational
complexity in the sampling.
For the GP model, look at the posterior given in
Equation~(\ref{eq-gpl-g}).  Equations~(\ref{eq-gpl2})
shows that the $\vec{l}_{(i)}$'s can be marginalised out.
Likewise, $\mat\Theta$ can be marginalised out because it is an instance of
a Dirichlet.
This yields a Gamma-Poisson posterior $p\left( \mat{V}_{(1)},\ldots,
\mat{V}_{(I)}\,| \,
\vec{\alpha},\vec{\beta},K,\mbox{GP}\right)$, with constants dropped:
\begin{equation}
\begin{split}
 \label{eq-lhgpf1} \prod_{i} \left( \prod_k
\frac{\Gamma(c_{k,(i)}+\alpha_k)}
                        {(1+\beta_k)^{c_{k,(i)}+\alpha_k}}
\prod_{j,k} \frac{1} {v_{j,k,(i)}! } \right)\\ \prod_k \frac{\prod_j
\Gamma\left(\gamma_j + \sum_i v_{j,k,(i)}\right) } {  \Gamma\left(
\sum_j \gamma_j + \sum_i c_{k,(i)}\right) }
\end{split}\end{equation}
Below it is shown that a short sampling routine can be based on
this.

A similar formula applies in the conditional GP case using
Equation~(\ref{eq-gpcl2}) for the marginalisation of $\vec{l}_{(i)}$'s.
The first term with $\prod_k$ in Equation~(\ref{eq-lhgpf1}) becomes
\[
 (1-\rho_k) \frac{\Gamma(c_{k,(i)}+\alpha_k)}{\Gamma(\alpha_k)}
              \frac{\beta_k^{\alpha_k}}
                        {(1+\beta_k)^{c_{k,(i)}+\alpha_k}}
 +\rho_k 1_{c_{k,(i)}=0} ~.
\]
Likewise a similar formula applies in the Dirichlet-multinomial
version using Equation~(\ref{eq-ddl2}):
\begin{equation}
\begin{split} \label{eq-lhgpf2} \prod_{i} \left( \prod_k
\Gamma(c_{k,(i)}+\alpha_k) \prod_{jk} \frac{1} {v_{j,k,(i)}! }
\right)\\ \prod_k \frac{\prod_j \Gamma\left(\gamma_j + \sum_i
v_{j,k,(i)}\right) } {  \Gamma\left( \sum_j \gamma_j + \sum_i
c_{k,(i)}\right) }
\end{split}\end{equation}
Here a term of the form $\Gamma\left( \sum_k(c_{k,(i)}+\alpha_k)
\right)$ drops out because $\sum_k c_{k,(i)}=L_{(i)}$ is known and
thus constant.

\comaWOff{If sequence of components \coma{sequence of
\underline{components}??? or words? that secion disccussed words
not components} formulation is used, instead of bag of components,
as discussed in Section~\ref{sct-bs}, then the $v_{j,k,(i)}!$
terms are dropped. \coma{I don't understand how would dropping
these terms affect anything in particular... If it's not
important, don't even mention it.}}

Now the posterior distributions have been marginalised
for each of the three models, GP, CGP and DM,
a Gibbs sampling scheme needs to be developed.
Each set $\{v_{j,k,(i)}\,:\,k \in 1,\ldots,K\}$ sums to
$w_{j,(i)}$, moreover the forms of the functions in
Equations~(\ref{eq-lhgpf1}) and~(\ref{eq-lhgpf2}) are
quite nasty.
A way out of this mess is to convert the scheme from
a bag of words model, implicit in the use of
$\mat{V}_{(i)}$ and $\vec{w}_{(i)}$, to a sequence of
words model.

This proceeds as follows.
Run along the $L_{(i)}$ words in a document and update
the corresponding component assignment for each word.
Component assignments for the $i$-th document
are in a $L_{(i)}$-dimensional
vector $\vec{k}_{(i)}$, where each entry takes a value
from $1,\ldots,K$.
Suppose the $l$-th word has word index $j_{l}$.
In one step, change the counts
$\{v_{j_l,k,(i)}\,:\,k \in 1,\ldots,K\}$ by one (one is increased
and one is decreased) keeping the total $w_{j_l,(i)}$ constant. For
instance, if a word is originally in component $k_1$ but
updating by Gibbs sampling to $k_2$,
then decrease $v_{j_l,k_1,(i)}$ by one and
increase $v_{j_l,k_2,(i)}$ by one. Do this for $L_{(i)}$ words
in the document, for each document.
Thus at word $l$ for the $i$-th document,
we sample component assignment $k_{l,(i)}$ according to the posterior
for $k_{l,(i)}$ with all other assignments fixed.
This posterior is proportional to
(the denominator is a convenient constant)
\[
\frac{\left.p\left( \mat{V}\,| \,\mbox{sequential},
\vec{\alpha},\vec{\beta},K,\mbox{GP}\right)
\right|_{v_{j_l,k,(i)}\leftarrow v_{j_l,k,(i)}+1_{k\neq k_l} }}
{\left.p\left( \mat{V}\,| \, \mbox{sequential},
\vec{\alpha},\vec{\beta},K,\mbox{GP}\right)
\right|_{v_{j_l,k_l,(i)}\leftarrow v_{j_l,k_l,(i)}-1}} ~,
\]
where the notation ``sequential'' is added to the right-hand side
because the combinatoric terms $v_{j,k,(i)}!$ of
Equation~(\ref{eq-lhgpf1}) need to be dropped.
This formula simplifies dramatically because $\Gamma(x+1)/\Gamma(x)=x$.

Derived sampling schemes are given in Table~\ref{tbl-rbgs}.  The
$(i)$ subscript is dropped and assumed for all counts, and $j=j_l$
is the word index for the word whose component index is being
resampled.  Since $k_l$ is being sampled, a $K$ dimensional
probability vector is needed.  The table gives the unnormalised
form.
\begin{table}
\begin{center}
\begin{tabular}{|l|p{3in}|}
\hline Model & ~~~~~~~~~~~~
Sampling Proportionality\\\hline GP &
\begin{minipage}{3in}
\vspace*{3pt} $$\frac{ \gamma_j + \sum_i v_{j,k} }
               { \sum_j \gamma_j + \sum_i c_{k} }
              \frac{c_{k}+\alpha_k}{1+\beta_k }$$
\end{minipage}\\[5pt]\hline
CGP &
\begin{minipage}{3in}
\vspace*{3pt} When $c_k>0$ use the proportionality of the
GP case, and otherwise
\[
\frac{ \gamma_j + \sum_i v_{j,k} }
               { \sum_j \gamma_j + \sum_i c_{k} }
\frac{\alpha_k}{1+\beta_k }
\frac{(1-\rho_k)\beta_k^{\alpha_k}}{(1-\rho_k)\beta_k^{\alpha_k} +
\rho_k(1+\beta_k)^{\alpha_k} }
\]
\end{minipage}\\[5pt]\hline
DM &
\begin{minipage}{3in}
\vspace*{3pt} $$\frac{ \gamma_j + \sum_i v_{j,k} }
               { \sum_j \gamma_j + \sum_i c_{k}
               }(c_{k}+\alpha_k).$$
\end{minipage}\\[5pt]\hline
\end{tabular}
\end{center}
\caption{Sampling $k_l=k$ given $j=j_l$ for Rao-Blackwellised Gibbs} \label{tbl-rbgs}
\end{table}

This Rao-Black\-well\-ised Gibbs algorithm is given in Figure~\ref{fig-rbga}.
\begin{figure*}[htb]
\begin{center}
\begin{boxedminipage}{5.8in}
\begin{enumerate}
\item Maintain the sufficient statistics for ${\mat\Theta}$, given
by $\sum_{i} v_{j,k,(i)}$ for each $j$ and $k$, and the sufficient
statistics for the component proportions
$\vec{l}_{(i)}$/$\vec{m}_{(i)}$ given by $\vec{c}_{(i)}$.
\item
For each document $i$, retrieve the $L_{(i)}$ component
 assignments for each word then:
\begin{enumerate}
\item Recompute statistics for $\vec{l}_{(i)}$/$\vec{m}_{(i)}$
 given by $c_{k,(i)} = \sum_j v_{j,k,(i)}$
for each $k$ from the individual component
assignment for each word.
\item For each word $l$ with word index $j_l$ and component
assignment $k_l$ in the document,  resample the component
assignment for this word according to the marginalised likelihoods
in this section.
\begin{enumerate}
\item First decrement $v_{j_l,k_l,(i)}$ and $c_{k_l,(i)}$ by one to remove
the component assignment for the word. \item Sample $k_l=k$
proportionally as in Table~\ref{tbl-rbgs}. \item Increment
$v_{j_l,k_l,(i)}$ and $c_{k_l,(i)}$.
\end{enumerate}
\item Concurrently, record the log-probability such as
$p\left( \vec{w}_{(i)}\,|
\,  \mat{V}_{(i)},
      \vec{\alpha},\vec{\beta},{\mat\Theta},K,\mbox{GP}\right)$
for the appropriate model.
 \item Concurrently,
 update the sufficient statistics for $\vec{l}_{(i)}$/$\vec{m}_{(i)}$
and ${\mat\Theta}$.
\end{enumerate}
\end{enumerate}
\end{boxedminipage}
\end{center}
\caption{One Major Cycle of Rao-Blackwellised Gibbs Algorithm for
DCA} \label{fig-rbga}
\end{figure*}
As before, an approximately unbiased log probability can be recorded
in Step~2(c). This requires a value for ${\mat\Theta}$. While the
sufficient statistics could be used to supply the current mean
estimate for ${\mat\Theta}$, this is not a true sampled quantity.
An alternative method is to make a sample of ${\mat\Theta}$ in
each major cycle and use this.

\subsubsection*{Implementation notes:}
Due to Rao-Black\-well\-isat\-ion, both the $\vec{l}_{(i)}$'s and
${\mat\Theta}$ are effectively re-estimated with each sampling
step, instead of once after the full pass over documents. This is
most effective during early stages, and explains the superiority
of the method observed in practice. Moreover, it means only one
storage slot for ${\mat\Theta}$ is needed (to store the sufficient
statistics), whereas in direct Gibbs two slots are needed (current
value plus the sufficient statistics). This represents a major
saving in memory. Finally, the $\vec{l}_{(i)}$'s and ${\mat\Theta}$ can be
sampled at any stage of this process (because their sufficient
statistics make up the totals appearing in the formula), thus
Gibbs estimates for them can be made as well during the MCMC
process.

\subsection{Historical notes}

Some previous algorithms can now be placed into context.
\begin{description}
\item[NMF: ]
Straight maximum likelihood, e.g. in \cite{lee99learning},
expressed in terms of Kullback-Leibler divergence
minimization, where optimisation jointly applies to the
latent variables
(see Section~\ref{sct-nmf}).
\item[PLSI: ] Annealed maximum likelihood \cite{hofmann99probabilistic}, best
viewed in terms of its clustering precursor such as by
\cite{hofmann97pairwise},
\item[Various Gibbs:]
Gibbs sampling on $\mat{V}_{(i)}$, $\vec{l}_{(i)}$/$\vec{m}_{(i)}$ and
 $\vec{\Theta}$ in turn using a full probability distribution
 by \cite{pritchard:1999}, or Gibbs sampling on $\mat{V}_{(i)}$ alone (or
 equivalently, component assignments for words in the sequence
 of words representation) after marginalising out
$\vec{l}_{(i)}$/$\vec{m}_{(i)}$
 and $\vec{\Theta}$ by \cite{griffiths04pnas},
\item[LDA: ]
variational approximation with Kullback-Leibler divergence
 by \cite{Blei03},
a significant introduction because of its speed.
\end{description}
Expectation propagation \cite{minka-uai02} requires $O(KS)$ latent
variables stored, a prohibitive expense compared to the $O(S)$ or
$O(KI)$ of other algorithms. Thus it has not been covered here.

\subsection{Other Aspects for Estimation and Use}

A number of other algorithms are needed to put these models into
regular use.

\paragraph{Component parameters: }
The treatment so far has assumed the parameter vectors $\vec\alpha$ and
$\vec\beta$ are given. It is more usual to estimate these
parameters with the rest of the estimation tasks as done by
\cite{Blei03,canny04}. This is feasible because the parameters are
shared across all the data, unlike the component vectors
themselves.

\paragraph{Estimating the number of components $K$:}
The number of components $K$ is usually a constant assumed a
priori. But it may be helpful to treat as a fittable parameter or
a random variable that adapts to the data.
In popular terms, this could be used to find
the ``right'' number of components, though in practice and theory
such a thing might not exist.
To obtain best-fitting
$K$, we can employ cross-validation, or we assess the
\emph{evidence} (or marginal likelihood) for the model given a
particular choice of $K$ \cite{Chib95,BuntineJakulin04}. In
particular, evidence is the posterior probability of the data
given the choice of $K$ after all other parameters have been
integrated out.

\paragraph{Use on new data:}
A typical use of the model requires performing inference related
to a particular document.  Suppose, for instance, one wished to
estimate how well a snippet of text, a query, matches a document.
Our document's components are summarised by the latent variables
$\vec{m}$ (or $\vec{l}$).  If the new query is represented by
$\vec{q}$, then $p(\vec{q}|\vec{m},{\mat\Theta},K,\mbox{GP})$ is the
matching quantity one  would like ideally. Since $\vec{m}$ is
unknown, we must average over it. Various methods have been
proposed \cite{minka-uai02,BuntineJakulin04}.

\comaOff{... We can employ the reduced description of the document in
terms of the component vector $\vec{c}$ instead of the description
in terms of the words..."}

\paragraph{Alternative components:}
Hierarchical components have been suggested
\cite{BuntineJakulin04} as a way of organising an otherwise large
flat component space.  For instance, the Wikipedia with over half
a million documents can easily support the discovery of several
hundred components. Dirichlet processes have been developed as an
alternative to the $K$-dimensional component priors in the
Dirichlet-multinomial/discrete model \cite{trespAIS2005}, although
in implementation the effect is to use $K$-dimensional Dirichlets
for a large $K$ and delete low performing components.

\comaOff{Uhm, dirichlet processes are really something that should be
mentioned above when discussing optimal $K$ - namely, through
Dirichlet processes we ``integrate out'' the choice of $K$ -
essentially performing `proper' Bayesian inference rather than
some sort of MAP determination of $K$. With Dirichlet processes
multiple choices of $K$ coexist in the same model. --- But a major
purpose of component analysis is to obtain interpretable
components, and it is unclear how to derive interpretations from
models with a probability distribution over $K$.}

\section{Applications}
\label{sct-app}

This section briefly discusses two applications of the methods.

\subsection{Voting Data}

One type of political science data are the \emph{roll calls}.
There were 459 roll calls in the US Senate in the year 2003. For
each of those, the vote of every senator was recorded in three
ways: `Yea', `Nay' and `Not Voting'. The outcome of the roll call
can be positive (e.g., Bill Passed, Nomination Confirmed)
corresponding to `Yea', or negative (e.g., Resolution Rejected,
Veto Sustained). Hence, the outcome of the vote can be interpreted
as the 101st senator, by associating positive outcomes with `Yea'
and negative outcomes with `Nay'.

\subsubsection{Application of the Method:}
We can now map the roll call data to the DCA framework. For each
senator $X$ we form two `words', where $w_{X,y}$ implies that $X$
voted `Yea', and $w_{X,n}$ implies that $X$ voted `Nay'. Each roll
call can be interpreted as a document containing a single
occurrence of some of the available words. The pair of words
$w_{X,y},w_{X,n}$ is then treated as a binomial, so the
multivariate formulation of Section~\ref{sct-mv} is used. Priors
for ${\mat\Theta}$ were Jeffreys priors, $\vec\alpha$ was
(0.1,0.1,...,0.1), and regular Gibbs sampling was used.

Special-purpose models are normally used for interpreting roll
call data in political science, and they often postulate a model
of rational decision making. Each senator is modelled as a
position or an \emph{ideal point} in a continuous spatial model of
preferences \cite{Clinton04}. For example, the first dimension
often delineates the liberal-conservative preference, and the
second region or social issues preference. The proximities between
ideal points `explain' the positive correlations between the
senators' votes. The ideal points for each senator can be obtained
either by optimization, for instance, with the optimal
classification algorithm~\cite{Poole00}, or through Bayesian
modelling~\cite{Clinton04}.

Unlike the spatial models, the DCA interprets the correlations
between votes through membership of the senators in similar blocs.
Blocs correspond to latent component variables. Of course, we can
speak only of the probability that a particular senator is a
member of a particular bloc. The corresponding probability vector
is normalized and thus assures that a senator is always a member
of one bloc on the average. The outcome of the vote is also a
member of several blocs, and we can interpret the membership as a
measure of how influential a particular bloc is.

Our latent senator (bloc) can be seen as casting votes in each
roll call. We model the behavior of such latent blocs across the
roll calls, and record it: it has a behavior of its own. In turn,
we also model the membership of each senator to a particular bloc,
which is assumed to be constant across all the blocs.

A related family of approaches is based on modelling relations or
networks using blocks or groups. There, a roll call would be
described by one network, individual senators would be nodes in
that network, and a pair of nodes is connected if the two senators
agreed. Discrete latent variables try to explain the existence of
links between entities in terms of senators' membership to blocks,
e.g., \cite{Holland83,Snijders97}.

Several authors prefer the block-model approach to modelling roll
call data \cite{Wang05}. The membership of senators to the same
block along with a high probability for within-block agreements
will explain the agreements between senators. While a bloc can be
seen as having an opinion about each issue, a block does not (at
least not explicitly). The authors also extended this model to
`topics', where the membership of senator to a particular block
depends on the topic of the issue; namely, the agreement between
senators depends on what is being discussed. The topic is also
associated with the words that appear in the description of an
issue.

\subsubsection{Visualization:}
We can analyze two aspects of the DCA model as applied to the
roll call data: we can examine the membership of senators in
blocs, and we can examine the actions of blocs for individual
issues. The approach to visualization is very similar, as we are
visualizing a set of probability vectors. We can use the gray
scale to mirror the probabilities ranging from 0 (white) to 1
(black).

As yet, we have not mentioned the choice of $K$ - the number of
blocs. Although the number of blocs can be a nuisance variable,
such a model is distinctly more difficult to show than one for a
fixed $K$. We obtain the following negative logarithms to the base
2 of the model's likelihood for $K = 4,5,6,7,10$: 9448.6406,
9245.8770, 9283.1475, 9277.0723, 9346.6973. We see that $K = 5$ is
overwhelmingly selected over all others, with $K = 4$ being far
worse. This means that with our model, we best describe the roll
call votes with the existence of five blocs. Fewer blocs do not
capture the nuances as well, while more blocs would not yield
reliable probability estimates given such an amount of data.
Still, those models are also valid to some extent. It is just that
for a single visualization we pick the best individual one of
them.

\begin{figure*}
\begin{center}
\epsfxsize=6.5in\epsfbox{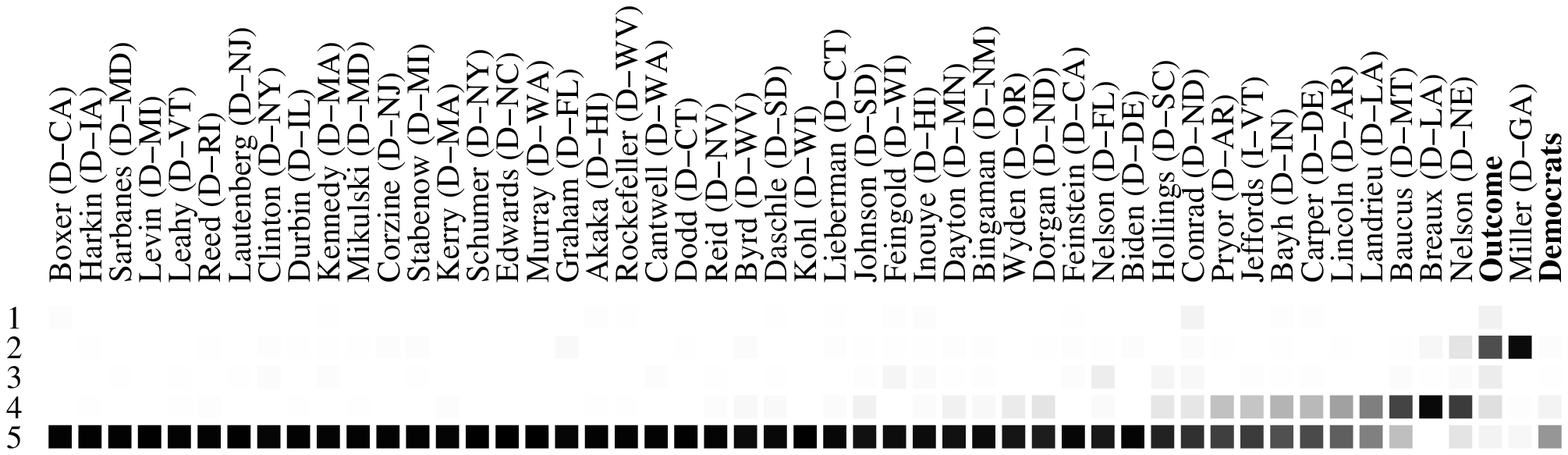}
\end{center} \vspace{-1.5cm} \caption{Component membership for Democrats}
\label{fig-dem}\vspace{5mm}
\begin{center}\epsfxsize=6.5in\epsfbox{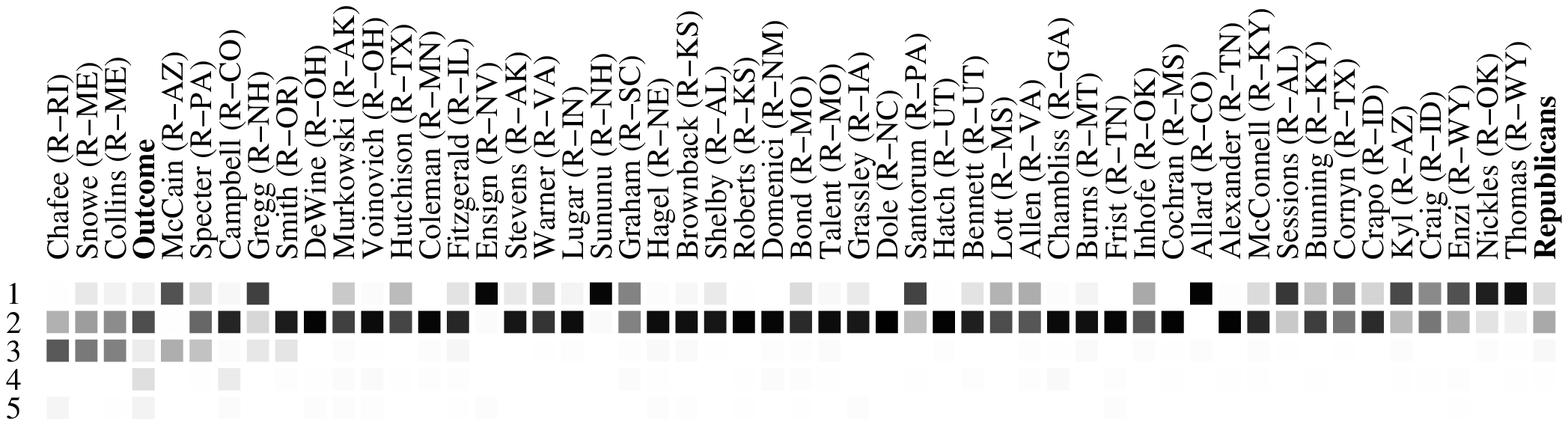}
\vspace{-1.5cm} \end{center} \caption{Component membership for
Republicans} \label{fig-rep}
\end{figure*}

We will now illustrate the membership of senators in blocs. Each
senator is represented with a vertical bar of 5 squares that
indicate his or her membership in blocs. We have arranged the
senators from left to right using the binary PCA approach of
\cite{deLeeuw03}. This ordering attempts to sort senators from the
most extreme to the most moderate and to the most extreme again.
Figure~\ref{fig-dem} shows the Democrat senators
and Figure~\ref{fig-rep} the Republicans.

We can observe that component 5 is the Democrat majority. It is
the strongest overall component, yet quite uninfluential about the
outcome. Component 4 are the moderate Democrats, and they seem
distinctly more influential than the Democrats of the majority.
Component 3 is a small group of Republican moderates. Component 2
is the Republican majority, the most influential bloc. Component 1
is the Republican minority, not very influential. Component 1
tends to be slightly more extreme than component 2 on the average,
but the two components clearly cannot be unambiguously sorted.

\section{Classification Experiments}

DCA is not trying to capture those aspects of the text that are
relevant for distinguishing one class of documents from another
one. Assume our classification task is to distinguish newspaper
articles on politics from all others. For DCA, modelling the
distribution of a word such as `the' is equally or more important than
modelling the distribution of a highly pertinent word such as
`tsunami' in classifying news reports. Of course, this is only a
single example of sub-optimality. Nevertheless, if there are several
methods of reducing the dimensionality of text that all disregard the
classification problem at hand, we can still compare them with respect
to classification performance.

We used MPCA and tested its use in its role as a feature
construction tool, a common use for PCA and ICA, and as a
classification tool. For this, we used the 20 newsgroups collection
described previously as well as the Reuters-21578 collection
\footnote{The Reuters-21578, Distribution 1.0 test collection is
available from David D. Lewis' professional home page, currently:
http://www.research.att.com/$\sim$lewis}. We employed the
$\textrm{SVM}^{light}$ V5.0 \cite{Joachims99} classifier with
default settings.  For classification, we added the class as a
distinct multinomial (cf. Section~\ref{sct-mv}) for the training
data and left it empty for the test data, and then predicted the
class value. Note that for performance and accuracy, SVM
\cite{Joachims98} is a clear winner \cite{Lewis04}. It is
interesting to see how MPCA compares.

Each component can be seen as generating a number of words in each
document. This number of component-generated words plays the same
role in classification as does the number of lexemes in the
document in ordinary classification. In both cases, we employed
the \tfidf\ transformed word and component-generated word counts as feature
values. Since SVM works with sparse data matrices, we assumed that
a component is not present in a document if the number of words
that a component would have generated is less than 0.01. The
components alone do not yield a classification performance that
would be competitive with SVM, as the label has no distinguished
role in the fitting. However, we may add these component-words in
the default bag of words, hoping that the conjunctions of words
inherent to each component will help improve the classification
performance.

For the Reuters collection, we used the ModApte split. For each of
the 6 most frequent categories, we performed binary classification.
Further results are disclosed in Table~2 \footnote{The numbers are
percentages, and `P/R' indicates precision/recall.}. No major change
was observed by adding 50 components to the original set of words.
By performing classification on components alone, the results were
inferior, even with a large number of components. In fact, with 300
components, the results were worse than with 200 components,
probably because of over-fitting. Therefore, regardless of the
number of components, the SVM performance with words cannot be
reproduced by component-generated words in this collection.
\begin{table}[h]
\caption{SVM Classification Results}\label{t:reut}
\begin{center}

\begin{tabular}{lrrrrr}\hline
 & \multicolumn{2}{c}{SVM} & & \multicolumn{2}{c}{SVM+MPCA}\\
 CAT & ACC. & P/R & &ACC. & P/R \\
 \hline

 earn & 98.58   & 98.5/97.1 && 98.45 & 98.2/97.1 \\
 acq & 95.54    & 97.2/81.9 && 95.60 & 97.2/82.2 \\
 moneyfx & 96.79& 79.2/55.3 && 96.73 & 77.5/55.9 \\
 grain & 98.94  & 94.5/81.2 && 98.70 & 95.7/74.5 \\
 crude & 97.91  & 89.0/72.5 && 97.82 & 88.7/70.9 \\
 trade & 98.24  & 79.2/68.1 && 98.36 & 81.0/69.8 \\ \hline \\[3mm]
\end{tabular}\\
\begin{tabular}{lrrrrr}\hline
 & \multicolumn{2}{c}{MPCA (50 comp.)}& & \multicolumn{2}{c}{MPCA (200
comp.)}\\
 CAT & ACC. & P/R & &ACC. & P/R \\
 \hline

 earn    & 96.94 & 96.1/94.6& & 97.06 & 96.3/94.8 \\
 acq     & 92.63 & 93.6/71.1& & 92.33 & 95.3/68.2 \\
 moneyfx & 95.48 & 67.0/33.0& & 96.61 & 76.0/54.7 \\
 grain   & 96.21 & 67.1/31.5& & 97.18 & 77.5/53.0  \\
 crude   & 96.57 & 81.1/52.4& & 96.79 & 86.1/52.4 \\
 trade   & 97.82 & 81.4/49.1& & 97.91 & 78.3/56.0 \\ \hline
\end{tabular}
 \end{center}
 \end{table}
Classifying newsgroup articles into 20 categories proved more
successful. We employed two replications of 5-fold cross
validation, and we achieved the classification accuracy of 90.7\%
with 50 additional MPCA components, and 87.1\% with SVM alone.
Comparing the two confusion matrices, the most frequent mistakes
caused by SVM+MPCA beyond those of SVM alone were predicting
talk.politics.misc as sci.crypt (26 errors) and talk.religion.misc
predicted as sci.electron (25 errors). On the other hand, the
components helped better identify alt.atheism and
talk.politics.misc, which were misclassified as talk.religion.misc
(259 fewer errors) earlier. Also, talk.politics.misc and
talk.religion.misc were not misclassified as talk.politics.gun (98
fewer errors). These 50 components were not very successful alone,
resulting in 18.5\% classification accuracy. By increasing the
number of components to 100 and 300, the classification accuracy
gradually increases to 25.0\% and 34.3\%. Therefore, many
components are needed for general-purpose classification.

From these experiments, we can conclude that components may help
with tightly coupled categories that require conjunctions of words
(20 newsgroups), but not with the keyword-identifiable categories
(Reuters). Judging from the ideas in \cite{Jakulin03b}, the
components help in two cases: a) when the co-appearance of two
words is more informative than the sum of informativeness of
individual appearances of either word, and b) when the appearance
of one word implies the appearance of another word, which does not
always appear in the document.

\section{Conclusion}

In this article, we have presented a unifying framework for
various approaches to discrete component analysis, presenting them
as a model closely related to ICA but suited for sparse discrete data.
We have shown the relationships
between existing approaches here such as NMF, PLSI, LDA, MPCA and
GaP. For
instance, NMF with normalised results corresponds to an approximate
 maximum likelihood method for LDA,
and GaP is the most general family of models.
We have also presented the different algorithms available for
three different cases, Gamma-Poisson, conditional Gamma-Poisson
(allowing sparse component scores), and Dirichlet-multinomial.
This extends a number of algorithms previous developed for
MPCA and LDA to the general Gamma-Poisson model.
Experiments with the {\sc Mpca} software\footnote{
http://www.componentanalysis.org} show that a typical 3GHz desktop
machine can build models in a few days
with $K$ in the hundreds for 3 gigabytes of text.

These models share many similarities with both PCA and ICA, and are
thus useful in a range of feature engineering tasks in machine
learning and pattern recognition. A rich literature is also emerging
extending the model in a variety of directions. This is as much
caused by the surprising performance of the algorithms, as it is by
the availability of general Gibbs sampling algorithms that allow
sophisticated modelling.

\comaOff{Maybe some additional discussion would be needed,}

\subsection*{Acknowledgments}
Wray Buntine's work was supported by the Academy of Finland under
the PROSE Project, by Finnish Technology Development Center (TEKES)
under the Search-Ina-Box Project, by the IST Programme of the
European Community under ALVIS Superpeer Semantic Search Engine
(IST-1-002068-IP) and the PASCAL Network of Excellence
(IST-2002-5006778). While performing the research described in this
paper, Aleks Jakulin was working at University of Ljubljana and at
Jo{\v z}ef Stefan Institute and was supported by the Slovenian
Research Agency and by the IST Programme of the European Community
under SEKT Semantically Enabled Knowledge Technologies
(IST-1-506826-IP). The {\sc Mpca} software used in the experiments
was co-developed by a number of authors, reported at the code
website. The experiments were supported by the document processing
environment and test collections at the CoSCo group, and the
information retrieval software {\sc Ydin} of Sami Perttu.

\bibliographystyle{alpha}
\bibliography{bo,extra,hiit,interactions}

\end{document}